\documentclass[10pt]{article}
\usepackage[top=4.2cm, bottom=4.2cm, left=4.0cm, right=4.0cm]{geometry}
\usepackage{amssymb}
\usepackage{latexsym}
\usepackage{parskip}
\usepackage{amsthm}
\usepackage[all]{xy}
\usepackage[T1]{fontenc}
\usepackage{amssymb}
\usepackage{latexsym}
\usepackage{amsthm}
\usepackage[all]{xy}
\usepackage[T1]{fontenc}

\theoremstyle{definition}
\newtheorem{remark}{Remark}[section]

\theoremstyle{theorem}
\newtheorem{thm}{Theorem} \newtheorem{prop}[remark]{Proposition}
\newtheorem{lemma}[remark]{Lemma} \newtheorem{cor}[remark]{Corollary}

\newtheorem{dfn}[remark]{Definition}

\newtheorem{lemmaapp}{Lemma}[section]
\newtheorem{thmapp}[lemmaapp]{Theorem}
\newtheorem{corapp}[lemmaapp]{Corollary}
\theoremstyle{definition}
\newtheorem{remarkapp}[lemmaapp]{Remark}

\newcommand{\rrvert}{\vert}
\newcommand{\rrVert}{\Vert}
\newcommand{\llvert}{\vert}
\newcommand{\llVert}{\Vert}

 \title{{\bf Quenched Invariance Principles for the Discrete Fourier Transforms of a Stationary Process\\}} \author{\large David Barrera\\{\vspace{0.2 cm}} \\ {Department of Mathematical Sciences} \\{  University of Cincinnati}\\ {\small PO Box 210025,
Cincinnati, Oh 45221-0025, USA.}\\{\small { } barrerjd@mail.uc.edu} \\{ }\\ {\small \&} \\{ } \\ Centre de Math\'{e}matiques Appliqu\'{e}es \\{ \'{E}cole Polytechnique} \\ {\small Route de Saclay 91128, Palaiseau Cedex, France.}\\ {\small {} david.barrera@polytechnique.edu}}

 \date{} 
 \begin{document} \maketitle

\begin{abstract}
In this paper, we study the asymptotic behavior of the normalized
{cadlag} functions generated by the discrete Fourier transforms of a
stationary centered square-integrable process, started at a point.

We prove that the quenched invariance principle holds for{ averaged}
frequencies under no assumption other than ergodicity, and that this
result holds also for almost every fixed frequency under a certain
generalization of the Hannan condition and a certain rotated form of the
Maxwell and Woodroofe condition which, under a condition of weak dependence that we specify, is guaranteed for a.e. frequency. If the process is in particular weakly mixing, our results describe the
asymptotic distributions of the normalized discrete Fourier transforms
at every frequency other than $0$ and $\pi$ under the generalized
Hannan condition.

We prove also that under a certain regularity hypothesis the conditional
centering is irrelevant for averaged frequencies, and that the same
holds for a given fixed frequency under the rotated Maxwell and
Woodroofe condition but not necessarily under the generalized Hannan
condition. In particular, this implies that the hypothesis of regularity is
not sufficient for functional convergence without random centering at
a.e. fixed frequency.

The proofs are based on martingale approximations and combine results
from Ergodic theory of recent and classical origin with approximation
results by contemporary authors and with some facts from Harmonic
Analysis and Functional Analysis.

\end{abstract}

\medskip

\textit{MSC2010 subject classification:} 60B10, 60F17, 60G10, 60G42, 37A05, 37A50.


\section{Introduction}
\label{int}

A recent result by Barrera and Peligrad (\cite{barpel}, Theorem~1)
shows that the quenched central limit theorem holds for the normalized components of the Fourier transforms of a stationary process in
$L^{2}$ orthogonal to the subspace of functions that are measurable with
respect to the initial sigma algebra. This quenched limit theorem
corresponds to previous, annealed limit theorems developed first by Wu
in \cite{Wu} and improved later by Peligrad and Wu in \cite{PeWu}.

The paper \cite{PeWu} contains also an invariance principle
(\cite{PeWu}, Proposition~2.1) for the Fourier Transforms. In
\cite{barpel}, on the other side, the problem of quenched convergence
for the sample paths generated by the discrete Fourier transforms is not
studied, which in particular leaves untouched the problem of extending
to the quenched setting the aforementioned invariance principle by Peligrad and Wu. 

In this paper, we will address this problem. We will prove an ``averaged
frequency'' quenched limit theorem (Theorem~\ref{Barinvpri})
corresponding to the invariance principle by Peligrad and Wu. We will
also see that the asymptotic behavior of the normalized sample paths
started at a point can be described for almost every fixed frequency
under a certain ``weak form'' of the Hannan condition (see
(\ref{imhancon})) and under a certain ``fixed frequency form'' of the
Maxwell and Woodroofe condition (see (\ref{maxwoo})).

More specifically, we will see that, under the weakened Hannan condition
(\ref{imhancon}), martingale approximations are possible for every
frequency other than (perhaps) zero, and we will deduce and describe
asymptotic distributions for the normalized sample paths at every
frequency other than those corresponding to the ``square root'' of the
point spectrum of the Koopman operator associated to the process. We
will also see that the same conclusion holds for a fixed frequency in
this set provided that (\ref{maxwoo}) holds, with the additional
conclusion that, in this case, the random centering is not needed. We will give a
sufficient condition, (\ref{equpromaxwooae}), for the fulfillment
of (\ref{maxwoo}) at a.e. frequency.

We emphasize that no assumption of regularity (see (\ref{reg})) is
needed for this, though some of the proofs are reached by first reducing
to the regular case. We also emphasize that our proofs take advantage
of some recently developed techniques in the realm of the calculus of asymptotic
distributions. See, for instance, Lemma~\ref{gengenbarinvpri} and its
proof, and the results by Cuny and Voln\'{y} in \cite{CuVo}. We point
out as well that the forthcoming results are valid for complex-valued
processes: the estimates needed for the martingale approximations do not
require the special properties (for instance, total ordering) of the
real numbers in an essential way,\footnote{The same observation applies
to many of the results referenced in this paper, which are often
stated under the unnecessary assumption that the processes under
consideration are real-valued.} and the martingale limit theorems are
valid as far as the square root of the point spectrum is avoided for the
rotations.\looseness=1

Our presentation is organized as follows: in Section~\ref{quecon}, we
will introduce some general facts and notions related to the quenched
convergence of stochastic processes. Then we will introduce, in Section~\ref{genset}, the setting in which our discussions will take place.
Section~\ref{thespaD} presents briefly the essentials of convergence in
distribution for {complex-valued} cadlag functions. In Section~\ref{resandcom}, we will present without proof our main results,
Theorems \ref{Barinvpri}, \ref{barinvprihancon} and \ref{promaxwooae},
preceded by a brief series of martingale approximation lemmas needed for
their proofs. Section~\ref{secmarcas} is devoted to the martingale case,
and Sections~\ref{prothe} and \ref{prothehancon} are devoted,
respectively, to the proofs of Theorems \ref{Barinvpri} and~\ref
{barinvprihancon} together with \ref{promaxwooae}. Finally, the {appendix} presents some results that are used along the previous proofs
and that deserve special mention due to their general or classical
nature.

\textit{Notation}. Throughout this paper, $\mathbb{N}=\{0,1,2,\dots\}$
denotes the natural numbers starting at zero, we also use the notation
$\mathbb{N}^{*}:=\mathbb{N}\setminus\{0\}$. Unless otherwise specified,
an expression of the form ``$\lim_{n}$'' (or ``$\limsup_{n}$'' or  ``$\liminf_{n}$'')
and ``$\rightarrow_{n}$'' must be read as ``$\lim_{n\to\infty}$'' (and
similarly for ``$\limsup_{n}$'' and ``$\liminf_{n}$'') and ``$
\rightarrow_{n\to\infty}$''.

\section{Quenched convergence}\label{quecon}
Let $(Y_{n})_{n}$ be a measurable sequence on some metric space
$(S,d)$. This is, for every $n$ ($n\in\mathbb{N}$ or $\mathbb{N}^{*}$
or $\mathbb{Z}$), $Y_{n}:(\Omega,\mathcal{F})\to(S,\mathcal{S})$ is
an $\mathcal{F}/\mathcal{S}$-measurable function where $(\Omega,
\mathcal{F})$ is a (fixed) measurable space and $\mathcal{S}$ is the Borel
sigma algebra of $S$. Let $\mathbb{P}$ be a given probability measure
on $(\Omega,\mathcal{F})$, so that $(\Omega,\mathcal{F},\mathbb{P})$
is a probability space, and denote by ``$\Rightarrow_{\mathbb{P}}$'' the
convergence in distribution with respect to $\mathbb{P}$.

The Portmanteau theorem (\cite{Bilconpromea}, Theorem~2.2) states,
among other equivalences, that if $Y:(\Omega',\mathcal{F}',
\mathbb{P}')\to(S,\mathcal{S})$ is a random element of $S$, then
$Y_{n}\Rightarrow_{\mathbb{P}} Y$ (as $n\to\infty$) if and only if for every continuous
and bounded function $f:S\to\mathbb{R}$
\begin{equation}
\label{usupro}
\int_{\Omega}f\circ Y_{n}(\omega)\, d\mathbb{P}(\omega )
\to_{n}
\int_{\Omega'}f\circ Y (z) \,d\mathbb{P}'(z),
\end{equation}
or, in the usual probabilistic notation, if and only if
%
\begin{equation}
\label{usupronot} \lim_{n}Ef(Y_{n})=Ef(Y),
\end{equation}
where $E$ is the expectation (Lebesgue integral) with respect to the
corresponding probability measures\footnote{When necessary, we will
indicate the underlying measures in some specified way, writing for
instance, ``$E_{\mathbb{P}}$ '' instead of ``$E$ '' for the
expectation with respect to $\mathbb{P}$.} and $f(Z):=f\circ Z$ (whenever this makes sense).

A~stronger kind of convergence, \textit{quenched} convergence, can be
defined in the following way: fix a sub-sigma algebra $\mathcal{F}
_{0}\subset\mathcal{F}$ representing, in a heuristic language, the
``initial information'' about (or the ``initial conditions of'') the
process $(Y_{n})_{n}$, and denote by $E_{0}$ the conditional expectation
with respect to $\mathcal{F}_{0}$. Then we will say that
\textit{$Y_{n}$ converges to $Y$ in the quenched sense with
respect to $\mathcal{F}_{0}$} if for every continuous and bounded
function $f:S\to\mathbb{R}$
%
\begin{equation}
\label{queconequconexp} E_{0} \bigl[f(Y_{n}) \bigr]
\to_{n} Ef(Y), \quad \quad\mathbb{P}\textit{-a.s.}
\end{equation}
Note that since this is pointwise convergence of
uniformly bounded functions (to a constant value), the dominated
convergence theorem guarantees that $\lim_{n}Ef(Y_{n})=Ef(Y)$, thus
quenched convergence implies convergence in distribution.

\medskip

\begin{remark}
\label{queimpsemque}
{The same argument, in combination with Theorem~34.2(v) in
\cite{BilProMea}, shows that if $\mathcal{G}_{0}\subset\mathcal{F}
_{0}$ is any sigma algebra, then the assumption of quenched convergence
with respect to $\mathcal{F}_{0}$ implies that $\lim_{n}E[f(Y_{n})|
\mathcal{G}_{0}]=Ef(Y)$, $\mathbb{P}$-a.s. In other words, \textit{quenched
convergence with respect to a given sigma algebra $\mathcal{F}_{0}$
implies quenched convergence with respect to any sub-sigma algebra of
$\mathcal{F}_{0}$}. Note also that one can interpret convergence in
distribution (or ``annealed'' convergence) as quenched convergence with respect to the trivial sigma algebra $\{\emptyset,\Omega\}$.}
\end{remark}

An example showing that the notion of quenched convergence is strictly
stronger than convergence in distribution can be constructed by starting
from any sequence $(Y_{n})_{n}$ of $\mathcal{F}_{0}$-measurable
functions and noticing that quenched convergence of $Y_{n}$ to $Y$ in
this case is the same as $f(Y_{n})\to Ef(Y)$, $\mathbb{P}$-a.s., for all
continuous and bounded functions $f$, which is not possible if, for
instance, $(Y_{n})_{n}$ takes the values $1$ and $0$ infinitely often,
$\mathbb{P}$-a.s.\footnote{Notice that, in this case $(f(Y_{n}))_{n}$
has no limit whatsoever, $\mathbb{P}$-a.s., for any $f$ with
$f(0)\neq f(1)$.}

More specifically, consider a sequence $(Y_{n})_{n}$ that converges in
distribution but gives $\mathbb{P}$-a.s. a sequence with infinitely many
$0$s and $1$s, and define $\mathcal{F}_{0}:=\sigma(\{Y_{n}\}_{n})$. For
instance, take the unit interval $[0,1]$ with the Lebesgue measure on
its Borel sigma algebra as the underlying probability space and let
$Y_{n}:[0,1]\to\{0,1\}$ be the characteristic function of $[0,1/2)$ or
the characteristic function of $[1/2,1]$ according to whether $n$ is,
respectively, even or odd. For another example, closely related to the content of this paper, the
reader is referred to \cite{baranexa}.

\medskip

Now assume that {$E_{0}$ is a \textit{regular conditional
expectation}: there exists a family of probability measures
$\{\mathbb{P}_{\omega}\}_{\omega\in\Omega}$ such that for every
integrable $X:(\Omega,\mathcal{F},\mathbb{P})\to\mathbb{R}$
%
\begin{equation}
\label{defe0reg} \omega\mapsto
\int_{\Omega}X(z)\,d\mathbb{P}_{\omega}(z)
\end{equation}
defines an $\mathcal{F}_{0}$-measurable version of $E_{0}X$.\footnote{More precisely we require, for a fixed version of $X$, the
existence of an $\mathcal{F}_{0}$-measurable set $\Omega_{X}$ with
$\mathbb{P}\Omega_{X}=1$ such that (\ref{defe0reg}) makes sense for
every $\omega\in\Omega_{X}$, and such that the function given by
(\ref{defe0reg}) if $\omega\in\Omega_{X}$, and zero otherwise, defines
an $\mathcal{F}_{0}$-measurable version of $E_{0}X$, i.e., an
$\mathcal{F}_{0}$-measurable function $\tilde{X}$ satisfying
\[
\int_{A}\tilde{X}(\omega)\,d\mathbb{P}(\omega)=
\int_{A}{X}(\omega)\,d \mathbb{P}(\omega)
\]
for every $A\in\mathcal{F}_{0}$. It is possible to prove the existence
of such $\Omega_{X}$ just by requiring that the function defined by
(\ref{defe0reg}) if $X\in L^{1}_{\mathbb{P}_{\omega}}$ and by zero
otherwise defines an $\mathcal{F}_{0}$-measurable version of
$E_{0}X$ for every $X\in L^{1}_{\mathbb{P}}$. See, for instance,
\cite{bardis}, Remark~11.2. The approximation argument for this
statement, which is similar to the one used to prove Lemma~\ref{mealem} in the {Appendix}, actually shows that it is enough to
require that for every $A\in\mathcal{F}$, $\omega\mapsto\mathbb{P}
_{\omega}A$ defines an $\mathcal{F}_{0}$-measurable version of
$\mathbb{P}(A|\mathcal{F}_{0})$.} The existence of such a family is
guaranteed if, for instance, $(\Omega,\mathcal{F})$ is a \textit{Borel space}, regardless of what $\mathcal{F}_{0}$ is. See Theorem~5.14 in
\cite{einwar}.

From now on, we will just say that $Y_{n}$ \textit{converges
to $Y$ in the quenched sense} to mean that the quenched convergence is
with respect to a fixed sigma algebra $\mathcal{F}_{0}$, returning to
the full description only if necessary to avoid ambiguity. Our first
result on quenched convergence is the following:

\medskip

\begin{prop}
\label{unicon}
Assume that $S$ is separable. If $E_{0}$ is regular and $Y_{n}$
converges to $Y$ in the quenched sense ($Y$ is defined on some
probability space $(\Omega',\mathcal{F}',\mathbb{P}')$), there exists
a set $\Omega_{0}\subset\Omega$ with $\mathbb{P}\Omega_{0}=1$ such
that for all $f:S\to\mathbb{R}$ continuous and bounded and all
$\omega\in\Omega_{0}$
%
\begin{equation}
\label{defunicon}
\int_{\Omega}f\circ Y_{n}(z)\,d\mathbb{P}_{\omega}(z)
\to_{n}
\int_{\Omega'}f\circ Y(z)\,d\mathbb{P}'(z).
\end{equation}
In particular, $Y_{n}$ converges to $Y$ in the quenched sense if and
only if for $\mathbb{P}$-a.e. $\omega$, $Y_{n}
\Rightarrow_{\mathbb{P}_{\omega}} Y$ as $n\to\infty$.
\end{prop}

Thus, we can choose the set of a.s. convergence in the definition of
quenched convergence \textit{uniform} over $\mathbf{C}^{b}(S)$: the space
of bounded, continuous functions $S\to\mathbb{R}$. The set
$\Omega_{0}$ depends, nonetheless, on $(Y_{n})_{n}$.
 
 {\bf Proof of Proposition \ref{unicon}:} Denote by $E^{\omega}$ the integration with respect to
$\mathbb{P}_{\omega}$, and consider functions
$U_{k,\epsilon}$ as in the statement~2. of Lemma~\ref{imppro} in the
{Appendix}. By the definition of regularity there exists, for all
$k\in\mathbb{N}$ and all $\epsilon\in\mathbb{Q}\cap(0,\infty)$, a
set $\Omega_{k,\epsilon}\subset\Omega$ with $\mathbb{P}
\Omega_{k,\epsilon}=1$ such that
\[
\Omega_{k,\epsilon}\subset\bigl\{\omega\in\Omega:E^{\omega}U_{k,
\epsilon}(Y_{n})
\to_{n} EU_{k,\epsilon}(Y)\bigr\}.
\]
Now take $\Omega_{0}:=\bigcap_{k,\epsilon} \Omega_{k,\epsilon}$ and
use Lemma~\ref{imppro}.
}.\qed
 
 The importance of Proposition~\ref{unicon} resides, for us, in the
following fact: to prove results on quenched convergence we will apply
some classical theorems to the processes involved in our arguments seen
as stochastic processes under the measures $\mathbb{P}_{\omega}$.
Without this result the uniformity of $\Omega_{0}$, which is eventually
necessary, would require a case-by-case approach, making the proofs much
longer and less transparent.

\medskip
\begin{remark}
\label{remproone}
We also point out the following: a set has $A\in\mathcal{F}$ satisfies
$\mathbb{P}A=1$ if and only if $\mathbb{P}_{\omega}A=1$ for
$\mathbb{P}$-a.e. $\omega$. To see this, use the equality
\[
\mathbb{P}A=
\int_{\Omega}\mathbb{P}_{\omega}A\,d\mathbb{P}(\omega ).
\]
\end{remark}

\section{General setting}\label{genset}

\subsection{Assumptions}\label{secass}
Our general setting, fixed from now on along this paper, will be the
following: first, denote by $\lambda$ the \textit{normalized Lebesgue measure} in the Borel sigma algebra $\mathcal{B}$
of $[0,2\pi)$. This is,
%
\begin{equation}
\label{norlebmeaequ} \lambda(A)=\frac{1}{2\pi} L(A)
\end{equation}
for every $A\in\mathcal{B}$, where $L$ is the Lebesgue measure on
$\mathcal{B}$.

By a \textit{random variable} we mean a \textit{complex-valued} measurable function $Y:\Omega\to\mathbb{C}$ defined
on some probability space $(\Omega,\mathcal{F},\mathbb{P})$.

Next, let $(X_{k})_{k\in\mathbb{Z}}$ be a strictly stationary, ergodic
sequence of random variables defined on a probability space
$(\Omega, \mathcal{F}, \mathbb{P})$. This is: { $X_{k}=X_{0}\circ T
^{k}$, where $T:\Omega\to\Omega$ is an ergodic, invertible, and
bimeasurable transformation}.

We will assume that $X_{0}\in L^{2}_{\mathbb{P}}(\mathcal{F}_{0})$ where
$\mathcal{F}_{0}\subset\mathcal{F}$ is a sigma algebra satisfying
$\mathcal{F}_{0}\subset T^{-1}\mathcal{F}_{0}$ (i.e. $T^{-1}$ is
$\mathcal{F}_{0}$-measurable), and we define $\mathcal{F}_{n}:=T^{-n}
\mathcal{F}_{0}$ for all $n\in\mathbb{Z}$, $\mathcal{F}_{-\infty}:=
\bigcap_{n\in\mathbb{Z}}\mathcal{F}_{n}$, and $\mathcal{F}_{\infty}:=
\sigma(\bigcup_{n\in\mathbb{Z}}\mathcal{F}_{n})$. Thus $(\mathcal{F}
_{k})_{k=-\infty}^{\infty}$ is an \textit{increasing
$T$-filtration}: $\mathcal{F}_{k}\subset\mathcal{F}_{k+1}$ and
$T^{-l}\mathcal{F}_{k}=\mathcal{F}_{k+l}$.

For any $n\in\{-\infty\}\cup\mathbb{Z}$, denote by $E_{n}$ the
conditional expectation with respect to $\mathcal{F}_{n}$, thus
$E_{n}Z:=E[Z|\mathcal{F}_{n}]$ for every $\mathbb{P}$-integrable $Z$,
and let the projection $\mathcal{P}_{n}$ be given by
\[\mathcal{P}_{n}Y:=E
_{n}Y-E_{n-1}Y.\]
Note that for $Y\in L^{2}_{\mathbb{P}}$, $\mathcal{P}
_{n}Y\in L^{2}_{\mathbb{P}}(\mathcal{F}_{n})\ominus L^{2}_{\mathbb{P}}(
\mathcal{F}_{n-1})$. This is, $\mathcal{P}_{n}Y$ is $\mathcal{F}_{n}$-measurable and $E_{n-1}\mathcal{P}_{n}Y=0$.

Assume also that{ $E_{0}$ is a regular conditional expectation}: as
explained in Section~\ref{quecon}, there exists a family of probability
measures $\{\mathbb{P}_{\omega}\}_{\omega\in\Omega}$ such that for
every integrable function $X$,
%
\begin{equation}
\label{conexpintreg} \omega\mapsto
\int_{\Omega}X(z)\,d\mathbb{P}_{\omega}(z)
\end{equation}
defines an $\mathcal{F}_{0}$-measurable version of $E_{0}X$ (see also
the footnote following (\pageref{defe0reg})).

We will finally assume that $\mathcal{F}_{0}$ and $\mathcal{F}$ are countably
generated or, alternatively, that $\mathcal{F}_{0}$ is countably
generated and $\mathcal{F}=\mathcal{F}_{\infty}$ (so $\mathcal{F}$ is \textit{also} countably generated). This alternative is
possible because we will deal only with
$\mathcal{F}_{\infty}$-measurable functions.


\medskip 

\begin{remark}
\label{remgenset}
If we also denote by $T:L^{1}_{\mathbb{P}}\to L^{1}_{\mathbb{P}}$ the
\textit{Koopman operator} associated to $T$, namely
$TY:=Y\circ T$ then, clearly, $X_{n}=T^{n} X_{0}$ for all $n$, and it
is not hard to see, using stationarity, that
%
\begin{equation}
\label{intkooopeequint} T^{r}E_{s}=E_{s+r}T^{r}
\end{equation}
(as operators in $L^{1}_{\mathbb{P}}$) for all integers $r$, $s$. Similarly, an application of the reverse
martingale convergence theorem (see for instance Theorem~5.6.1 and
Exercise 5.6.1 in \cite{dur}) shows that, for every $n\in
\mathbb{Z}$,
%
\begin{equation}
\label{intkooopeequinf} T^{n}E_{-\infty}=E_{-\infty}T^{n}.
\end{equation}
It is important to point out also that, again by the reverse martingale
convergence theorem, the following holds: for every $X\in L^{p}_{
\mathbb{P}}$ ($p\geq1$)
%
\begin{equation}
\label{revmarcon} E_{-\infty}X=\lim_{n} E_{-n}X,
\end{equation}
$\mathbb{P}$-a.s. and in $L^{p}_{\mathbb{P}}$.
\end{remark}

\medskip

\begin{remark}
\label{remkooope}
We also recall the following fact about the Koopman operator $T$: under
ergodicity, the eigenvalues of $T$ form a subgroup of $\mathbb{T}$, the
unit circle seen as a (Lie) group under the operation of multiplication
of complex numbers (see \cite{eisfarhaanag}, Proposition~7.17).
We will denote this group by $\mathrm{Spec}_{p}(T)$, the \textit{point spectrum} of $T$. Note that since
$L^{2}_{\mathbb{P}}$ admits a countable orthonormal basis ($
\mathcal{F}$ is countably generated) and the eigenspaces of $T$ are
mutually orthogonal ($T$ is measure preserving), $\mathrm{Spec}_{p}(T)$ is
countable. In particular,
%
\begin{equation}
\label{lammeapoispezer} \lambda \bigl( \bigl\{\theta\in[0,2\pi ): e^{i\theta}
\in \mathrm{Spec}_{p}(T) \bigr\}\bigr)=0.
\end{equation}
\end{remark}

\subsection{Quenched convergence in the product space}\label{queconprospa}
There is a special form of quenched convergence that will be of interest
to us: let ${\mathcal{G}}_{0}:=\mathcal{B}_{0}\otimes\mathcal{F}_{0}
\subset\mathcal{B}\otimes\mathcal{F}$ where $\mathcal{B}_{0}\subset
\mathcal{B}$ is a given sigma algebra (we will choose $\mathcal{B}
_{0}=\mathcal{B}$ or $\mathcal{B}_{0}=\{\emptyset,[0,2\pi)\}$
according to the problem under consideration). Assuming that $E[\,\cdot\,|\mathcal{B}_{0}]$ (conditional expectation
with respect to $\lambda$) is regular with regular measures
$\{\lambda_{\theta}\}_{\theta\in[0,2\pi)}$, an application of
Corollary~\ref{cormealem} in the {appendix} shows that, for any
$\mathcal{B}\otimes\mathcal{F}$-integrable function $Y=Y(\theta,
\omega)$, a version of $E[Y|\mathcal{G}_{0}]$ is given by
\[
E[Y|\mathcal{G}_{0}](\theta,\omega)=
\int_{[0,2\pi)}E^{\omega}\bigl[Y(x, \cdot)\bigr]\,d
\lambda_{\theta}(x),
\]
where $\omega\mapsto E^{\omega}Y(x,\cdot)$ is the version of
$E_{0}[Y(x,\cdot)]$ given by (\ref{conexpintreg}). 

From now on, we will
denote by $\tilde{E}_{0}$ the conditional expectation $E[\,\cdot\,|
\mathcal{G}_{0}]$. Whenever needed, we will work under the regular
version of $\tilde{E}_{0}$, so that for any integrable $Y=Y(\theta,
\omega)$
%
\begin{equation}
\label{regconexg0} \tilde{E}_{0}[Y](\theta,\omega)=
\int_{[0,2\pi)}
\int_{\Omega}Y(x,z) \,d\mathbb{P}_{\omega}(z)\,d
\lambda_{\theta}(x).
\end{equation}
%
\medskip

\begin{remark}
\label{remdouint}
Note that if $\mathcal{B}_{0}=\{\emptyset, [0,2\pi)\}$ is the trivial
sigma algebra, then $\lambda_{\theta}:=\lambda$ for all $\theta
\in[0,2\pi)$ defines a family of regular measures for $E[\,\cdot\,|
\mathcal{B}_{0}]$. And that in this case $\tilde{E}_{0}Y$, being constant in $\theta$ for
a fixed $\omega$, defines an $\mathcal{F}_{0}$-measurable function.
In particular, $\tilde{E}_{0}Y$ can be regarded, in this case, as a
random variable defined on $(\Omega,\mathcal{F},\mathbb{P})$.
\end{remark}

\section{The space $D[[0,\infty),\mathbb{C}]$}\label{thespaD}

This paper deals with convergence in distribution, under several
measures, of random elements of $D[[0,\infty),\mathbb{C}]$: the space
of functions $f:[0,\infty)\to\mathbb{C}$ that are continuous from the
right and have left-hand limits at every point (\textit{complex-valued} cadlag \textit{functions}). This space is an algebra
with the operation of multiplication and addition given by the usual
pointwise operations between complex-valued functions, and it is a ($
\mathbb{C}$ or $\mathbb{R}$-)vector space with the usual operation of
multiplication by constants regarded as constant functions.

Let $D[0,\infty)$ denote the space of real-valued cadlag functions $f:[0,\infty)\to \mathbb{R}$ (as presented for instance in \cite{Bilconpromea}, Section~16).
If we denote by $( \mbox{Re}(f),  \mbox{Im}(f))$ the
vector of real and imaginary parts
of a function $f\in D[[0,\infty),\mathbb{C}]$, the bijection
$D[[0,\infty),\mathbb{C}]\to D[0,\infty)\times D[0,\infty)$ given by
\[
f\mapsto\bigl( \mbox{Re}(f),  \mbox{Im}(f)\bigr)
\]
allows us to regard $D[[0,\infty),\mathbb{C}]$ as a topological space
whose topology is the topology generated by the product Skorohod
topology of $D[0,\infty)$.\footnote{But remember that this is
\textit{not} a topological vector space. See exercise 12.2 in
\cite{Bilconpromea}.} 

In particular $D[[0,\infty),\mathbb{C}]$ is
separable and metrizable. If we use the product metric
\[
d(f,g):=\bigl(\bigl(d\bigl( \mbox{Re}(f),  \mbox{Re}(g)\bigr)
\bigr)^{2}+\bigl(d\bigl( \mbox{Im}(f), \mbox{Im}(g)\bigr)
\bigr)^{2}\bigr)^{1/2},
\]
where $d$ (at the right) is the Skorohod distance defined by
\cite{Bilconpromea}, (12.16) then $D[[0,\infty),\mathbb{C}]$ is
complete. A~similar construction shows the corresponding facts for the
space $D[[0,m],\mathbb{C}]$ of functions that are restrictions of
elements in $D[[0,\infty),\mathbb{C}]$ to the interval $[0,m]$.

Denote by $\mathcal{D}_{\infty,\mathbb{C}}$ the Borel sigma algebra of
$D[[0,\infty),\mathbb{C}]$. The following observations will suffice to
clarify the proofs of convergence in distribution given here:

\begin{enumerate}
\item[1.] First, if $Y:(\Omega,\mathcal{F},\mathbb{P})\to D[[0,\infty),
\mathbb{C}]$ is a random function (this is, an $\mathcal{F}/
\mathcal{D}_{\infty,\mathbb{C}}$ measurable function, where
$(\Omega,\mathcal{F},\mathbb{P})$ is a probability space), then the
inequalities
\[
\mathbb{P}[Y_{n}\notin K_{1}\times K_{2}]\leq
\mathbb{P}\bigl[ \mbox {Re}(Y_{n}) \notin K_{1}\bigr]+
\mathbb{P}\bigl[ \mbox{Im}(Y_{n})\notin K_{2}\bigr]
\leq 2\mathbb{P}[Y_{n} \notin K_{1}\times K_{2}]
\]
show that \textit{a sequence of random elements $(Y_{n})_{n}$ in
$D[[0,\infty),\mathbb{C}]$ is tight if and only if $(\mbox{\upshape Re}(Y_{n}))_{n}$
and $( \mbox{{\upshape Im}}(Y_{n}))_{n}$ are tight}.
\item[2.] By an adaptation of the arguments in \cite{Bilconpromea}, it is
possible to show that the finite dimensional distributions are a
\textit{separating class} in $D[[0,\infty),\mathbb{C}]$: if for every
$t$ we denote by $\pi_{t}$ the coordinate function $\pi_{t}(x):=x(t)$,
then two probability measures $\mathbb{P}_{1}$ and $\mathbb{P}_{2}$ in
$D[[0,\infty),\mathbb{C}]$ coincide if and only if there exists a dense
subset $J\subset[0,\infty)$ such that for every $0\leq t_{1}\leq
\cdots\leq t_{n}$ in $J$ the $n$th dimensional distributions
$(\mathbb{P}_{j}\pi_{t_{1}}^{-1},\ldots,\mathbb{P}_{j}\pi_{t_{n}}
^{-1})$ ($j=1,2$) on $\mathbb{C}^{n}=\mathbb{R}^{2n}$ are the same. In
particular, one can prove that $\mathbb{P}_{n}\Rightarrow_{n}\mathbb{P}$
(``$\Rightarrow$'' denotes weak convergence of measures) by proving
tightness and convergence of the finite-dimensional distributions, and
if $m> 0$ is such that
%
\begin{equation}
\label{conconpoi} \mathbb{P} \Bigl\{x:\lim_{t\to m^{-}}x(t)\neq x(m)
\Bigr\}=0
\end{equation}
and $r_{m}: D[[0,\infty),\mathbb{C}]\to D[[0,m],\mathbb{C}]$ is the
restriction operator\footnote{The space $D[[0,m],\mathbb{C}]$ is defined via the product metric as above
starting from the space $D[0,m]$ of real-valued cadlag functions $f:[0,m]\to\mathbb{R}$. See (again)
Section~16 in \cite{Bilconpromea} or Section~7 in \cite{bardis} for more details on the last space.}
($r_{m}x(t)=x(t)$), then $\mathbb{P}_{n}\Rightarrow
\mathbb{P}$ in $D[[0,\infty),\mathbb{C}]$ implies that $\mathbb{P}
_{n}r_{m}^{-1}\Rightarrow\mathbb{P}r_{m}^{-1}$ in $D[[0,m],
\mathbb{C}]$. Conversely, if $(m_{k})_{k}$ is a sequence increasing to
infinity such that (\ref{conconpoi}) holds for all $m=m_{k}$, then
$\mathbb{P}_{n}r_{m_{k}}^{-1}\Rightarrow\mathbb{P}r_{m_{k}}^{-1}$ (on
$D[[0,m_{k}],\mathbb{C}]$) for all $k$ implies that $\mathbb{P}_{n}
\Rightarrow_{n}\mathbb{P}$.
\end{enumerate}

Let us finally mention that, by our definition of the topology of
$D[[0,\infty),\mathbb{C}]$, to prove that a function $Y:\Omega
\to D[[0,\infty),\mathbb{C}]$ defined on a measurable space
$(\Omega,\mathcal{F})$ is $\mathcal{F}/\mathcal{D}_{\infty,
\mathbb{C}}$ measurable it suffices to see the measurability of the real
and imaginary parts of $Y$. This observation, combined with the argument
in \cite{Bilconpromea}, p.~84, and with Theorem~16.6 in that book shows
that $Y$ is $\mathcal{F}/\mathcal{D}_{\infty,\mathbb{C}}$-measurable
if for every $t\in[0,\infty)$, $\omega\mapsto Y(\omega)(t)$ is
$\mathcal{F}$-measurable.

\section{Results and comments}\label{resandcom}
In this section, we will present the main results of this
paper. Some of the proofs are not difficult and can be given after the
statements. More technically demanding facts will be deferred to later
sections.

Let us start by introducing the notion of discrete Fourier transforms.

\medskip

\begin{dfn}[Discrete Fourier Transforms of a Stationary Process]
\label{defdisfoutra}
Let $T:\Omega\to\Omega$ be an invertible, bimeasurable,
measure-preserving transformation on a probability space $(\Omega,
\mathcal{F},\mathbb{P})$, let $Y_{0}:\Omega\to \mathbb{C}$ be a random variable, and let $Y_{k}:=T^{k}Y_{0}$ ($k\in
\mathbb{N}$ or $\mathbb{Z}$). For every $\theta\in[0,2\pi)$ and
$n\in\mathbb{N}^{*}$, the {$n$th \it discrete Fourier Transform of
$(Y_{k})_{k}$ at the frequency $\theta$} is defined by
%
\begin{equation}
\label{nfoutradef} S_{n}(Y_{0}, T, \theta,\omega):=\sum
_{k=0}^{n-1}e^{ik\theta}Y_{k}( {
\omega}).
\end{equation}
When $Y_{0}$ and $T$ are fixed and $\theta\in(0,2\pi)$ is given, we
will denote by $S_{n}(\theta)$ (or $S_{n}(\theta,\cdot)$) the random
variable $S_{n}(Y_{0},T,\theta,\cdot)$. If in addition $\theta=0$,
we denote by $S_{n}$ the random variable $S_{n}(Y_{0},T,0,\cdot)$. So
$S_{n}(\omega):=\sum_{k=0}^{n-1}Y_{k}(\omega)$.
\end{dfn}

\medskip

\begin{remark}
We strengthen the fact that, at this point, $(\Omega,\mathcal{F},
\mathbb{P})$ and $T$ in this definition are \textit{not} necessarily the
objects fixed in Section~\ref{secass}, though we will not go ``too far
away'' from them,\footnote{Indeed, we will need to work with product spaces and product maps along
some of the proofs and this makes convenient to relax these objects in
Definition~\ref{defdisfoutra}.} hence the specification of $T$ in the notation.

Nonetheless, \textit{throughout the rest of this section,
we will work under the assumptions specified in Section}~\ref{secass}.
In particular, we will not specify $T$ in the notation for discrete
Fourier transforms.
\end{remark}

We will address the problem of quenched convergence for the {cadlag}
random functions generated by the discrete Fourier transforms of
$(X_{k})_{k}$. Let us start by recalling the following theorem, proved
in \cite{barpel}:

\medskip

\begin{thm}\label{BaPeTh}
There exists a set $I\subset(0,2\pi)$ with $\lambda(I)=1$ such that,
for all $\theta\in I$, the random variables
%
\begin{equation}
\label{norparsumc} \frac{1}{\sqrt{n}} \bigl(S_{n}(\theta
)-E_{0}S_{n}(\theta) \bigr)
\end{equation}
converge in the quenched sense, as $n\to\infty$, to a complex Gaussian
random variable with independent real and imaginary parts, each with
mean zero and variance $\sigma^{2}(\theta)/2$, where
%
\begin{equation}
\label{asyvar} \sigma^{2}({\theta})=\lim_{n}
\frac{E_{0}\llvert S_{n}(\theta)-E_{0}S_{n}(
\theta)\rrvert^{2}}{n}=\lim_{n} \frac{E|S_{n}(\theta
)-E_{0}S_{n}(\theta)|^{2}}{n},
\end{equation}
$\mathbb{P}$-a.s. (thus the limit exists with probability one, and it
is nonrandom).
\end{thm}

We will refer to this as the \textit{quenched central limit theorem for the
discrete Fourier transforms of a stationary process}. Our main goal in
this paper is to explore possible extensions of this theorem to
corresponding quenched invariance principles.

\medskip
\begin{remark}
\label{remonsigsquathe}
The function $\theta\mapsto\sigma^{2}(\theta)$ defines a version of
the spectral density, with respect to~$\lambda$, of the process
$(X_{k}-E_{-\infty}X_{k})_{k}$. For details on this, see the proof of
Theorem~15.1 in \cite{bardis} (or combine Remark~\ref{remasyvar} and
(\ref{redregcas}) below with the proof of Theorem~3 in
\cite{lifpel}).
\end{remark}

In order to discuss the issue of quenched functional convergence, let
us define the following random variables, which will appear again later
along the proofs:
%
\begin{equation}
\label{firmardifr} D_{r,0}(\theta):=\sum_{k=0}^{r}
\mathcal{P}_{0}X_{k}e^{ik\theta},
\end{equation}
and
%
\begin{equation}
\label{marr} M_{r,n}(\theta):=\sum_{k=0}^{n-1}T^{k}D_{r,0}(
\theta)e^{ik\theta}.
\end{equation}
Where $0\leq r < \infty$ and we allow the value $r=\infty$ when
$D_{\infty,0}(\theta)$ makes sense as the limit ($r\to\infty$)
$\mathbb{P}$-a.s. and in $L^{2}_{\mathbb{P}}$ of $D_{r,0}(\theta)$.
Denote also by $V_{r,n}(\theta)$ the random function
%
\begin{equation}
\label{appfunr} V_{r,n}(\theta) (t):=\frac{1}{\sqrt{n}}M_{r, \lfloor
nt \rfloor}(
\theta),
\end{equation}
and by $W_{n}(\theta)$ the random function
%
\begin{equation}
\label{funqueinvpri} W_{n}(\theta) (t):=\frac{S_{ \lfloor nt \rfloor
}(\theta)-E
_{0}S_{ \lfloor nt \rfloor}(\theta)}{\sqrt{n}}.
\end{equation}

We will also consider the \textit{non-centered} version of $W_{n}(\theta
)$
%
\begin{equation}
\label{funqueinvprinoncen} U_{n}(\theta) (t):=W_{n}(\theta) (t)+
\frac{E_{0}S_{ \lfloor nt \rfloor}(
\theta)}{\sqrt{n}}=\frac{S_{ \lfloor nt \rfloor}(\theta)}{
\sqrt{n}}.
\end{equation}

\begin{remark}
\label{remtwopar}
When necessary, \textit{especially when discussing quenched convergence in
the product space} $([0,2\pi)\times\Omega,\mathcal{B}\otimes
\mathcal{F})$, we will specify the dependence on $\theta\in[0,2
\pi)$ \textit{and} $\omega\in\Omega$ by seeing these random elements
as processes with two parameters (see Theorem~\ref{Barinvpri} for
example). So for instance,
\[
D_{r,0}(\theta,\omega)=\sum_{k=0}^{r}
\mathcal{P}_{0}X_{k}(\omega)e ^{ik\theta},
\]
and so on.
\end{remark}

\subsection{Approximation lemmas}\label{genrescon}
Recall the definition of the set $I$ specified in Theorem~\ref{BaPeTh}. See \cite{barpel}; p.~289, 290: $\theta\in I$ if and only
if
\begin{enumerate}\item[1.] $e^{-2i\theta}\notin \mathrm{Spec}_{p}(T)$ (equivalently,
$e^{2i\theta}
\notin \mathrm{Spec}_{p}(T)$),
\item[2.] $E[\sup_{r}|\mathcal{P}_{0}S_{r}(\theta)|^{2}]< \infty$, and
\item[3.] $D_{\infty,0}(\theta):=\lim_{n} \mathcal{P}_{0}S_{n}(\theta
)$ exists
$\mathbb{P}$-a.s.
\end{enumerate}
Notice that, for $\theta\in I$, $D_{\infty,0}(\theta)$ is (also)  the limit in $L^{2}_{\mathbb{P}}$ of $D_{r,0}(\theta)$: {first, $D_{\infty,0}(\theta)\in L^{2}_{\mathbb{P}}$ by the dominated convergence theorem and 2., and since 
$D_{r,0}(\theta)-D_{\infty,0}(\theta)\to 0$ $\mathbb{P}$-a.s and $|D_{r,0}(\theta)-D_{\infty,0}(\theta)|^{2}\leq 2(|\sup_{r}D_{r,0}(\theta)|^2+|D_{\infty,0}(\theta)|^2)$ the dominated convergence theorem again implies that $D_{r,0}(\theta)\to_{r} D_{\infty,0}(\theta)$ in $L^{2}_{\mathbb{P}}$.}

\medskip

Notice that, for $\theta\in I$, $D_{\infty,0}(\theta)$ is (also) the
limit in $L^{2}_{\mathbb{P}}$ of $D_{r,0}(\theta)$: {first,
$D_{\infty,0}(\theta)\in L^{2}_{\mathbb{P}}$ by the dominated
convergence theorem and 2., and since $D_{r,0}(\theta)-D_{\infty,0}(
\theta)\to0$, $\mathbb{P}$-a.s. and $|D_{r,0}(\theta)-D_{\infty,0}(
\theta)|^{2}\leq2(\sup_{r}|D_{r,0}(\theta)|^{2}+|D_{\infty,0}(
\theta)|^{2})$ the dominated convergence theorem again implies that
$D_{r,0}(\theta)\to_{r} D_{\infty,0}(\theta)$ in $L^{2}_{
\mathbb{P}}$.}

\medskip

\begin{lemma}
\label{convnr}
Let $ \mbox{\upshape Re}(B), \mbox{\upshape Im}(B):\Omega'\to D[[0,\infty
)]$ be independent standard
Brownian motions\footnote{The notation may be a little shocking. It
is inspired in the fact that $B:= \mbox{Re}(B)+i \mbox{Im}(B)$ is a
$2$-dimensional
Brownian motion.} defined on some probability space $(\Omega',
\mathcal{F}',\mathbb{P}')$, and consider the set of frequencies
$\theta$ specified by condition 1. above. Then for all such
$\theta$ and all $0\leq r<\infty$, $V_{r,n}(\theta)$ given by
(\ref{appfunr}) converges in the quenched sense, as $n\to\infty$, to
%
\begin{equation}
\label{limindbror} B_{r}(\theta):= \bigl(E \bigl\llvert
D_{r,0}( \theta) \bigr\rrvert^{2}/2 \bigr)^{1/2}
\bigl(  \mbox{\upshape Re}(B)+i\, \mbox{\upshape Im}(B) \bigr).
\end{equation}
If $\theta\in I$ ($I$ is the set specified by 1.--3. above), then
$V_{\infty,n}(\theta)$ converges in the quenched sense to\vspace*{-1pt}
%
\begin{equation}
\label{limindbro} B(\theta):= \bigl(E \bigl\llvert D_{\infty,0}(\theta ) \bigr
\rrvert^{2}/2 \bigr)^{1/2} \bigl( \mbox{\upshape Re}(B)+i\,
\mbox {\upshape Im}(B) \bigr).
\end{equation}
\end{lemma}

{\bf Proof:} This is an immediate consequence of Theorem~\ref{marcas}
in Section~\ref{secmarcas}.
\qed

\medskip

Our results will make use as well of the following ``hypothetical''
quenched invariance principle:

\medskip
\begin{lemma}
\label{gengenbarinvpri}
If for a given $\theta\in[0,2\pi)$, (\ref{appfunr}) converges in the
quenched sense (with respect to $\mathcal{F}_{0}$) as $n\to\infty$ to
a random function $B_{r}(\theta)$ for all but finitely many $r$'s and\vspace*{-1pt}
%
\begin{equation}
\label{appiner} \lim_{r}\limsup_{n}E_{0}
\biggl[\frac{1}{n}\max_{1\leq k\leq n} \bigl\llvert
S_{k}( \theta)-E_{0}S_{k}(\theta)-M_{r,k}(
\theta) \bigr\rrvert^{2} \biggr]=0
\end{equation}
$\mathbb{P}$-a.s., then there exists a random function $B(\theta)$ such
that $B_{r}(\theta)$ converges in distribution to $B(\theta)$ as
$r\to\infty$, and (\ref{funqueinvpri}) converges to $B(\theta)$ as
$n\to\infty$ in the quenched sense.

If, in particular, $e^{2i\theta}\notin \mathrm{Spec}_{p}(T)$, then there exists
a nonnegative number $\sigma(\theta)$ characterized by\vspace*{-1pt}
%
\begin{equation}
\label{sigthe} \sigma^{2}(\theta):=\lim_{r}{E
\bigl\llvert D_{r,0}(\theta) \bigr\rrvert^{2}},
\end{equation}
and\vspace*{-1pt}
%
\begin{equation}
\label{auxdefbinf} B(\theta)=\frac{\sigma(\theta)}{\sqrt{2}} \bigl(\mbox{Re}(B)+i\,\mbox{Im}(B)
\bigr)
\end{equation}
(in distribution).
\end{lemma}

Note that if one agrees to call ``standard'' a $d$-dimensional Brownian
motion $B=(B_{1},\dots,B_{d})$ for which the entries are (independent
and) identically distributed (centered) Brownian motions and such that
$E|B(1)|^{2}=1$ where $|\cdot|$ is the Euclidean norm, the limit
function (\ref{auxdefbinf}) is just a standard $2$-dimensional Brownian
motion rescaled by $\sigma(\theta)$.

\medskip
\begin{remark}
If we assume \textit{a priori} the existence of $B(\theta)$ with
$B_{r}(\theta)\Rightarrow B(\theta)$ as $r\to\infty$, then the parts
of the proof of this and analogous results that rely in Theorem~\ref{limlimlem} in the {appendix} can still be carried over using instead
\cite{Bilconpromea}, Theorem~3.2.

This is possible for the values of $\theta$ that we will treat along
the proofs of the results in this paper: we will be restricted to cases
in which $D_{\infty,0}(\theta)$ is well defined, but since our
statement covers an (at least formally) more general case and the proofs
via the result in \cite{Bilconpromea} would be just the same (with the
estimates used here), we proceed via Theorem~\ref{limlimlem}.
\end{remark}

\medskip

\begin{remark}
\label{remasyvar}
It is important to remark also the following: in the context of Lemma~\ref{gengenbarinvpri}, (\ref{appiner}) implies that, necessarily,
%
\begin{equation}
\label{sigsquasyvar} \sigma^{2}(\theta)=\lim_{n}E_{0}
\frac{1}{n} \bigl\llvert S_{n}(\theta)-E_{0}S
_{n}(\theta) \bigr\rrvert^{2}
\end{equation}
$\mathbb{P}$-a.s.

Indeed, Lemma~\ref{dedmerpel} in the {appendix} implies, by orthogonality,
that
%
\begin{equation}
\label{expdrasanave} E \bigl[ \bigl\llvert D_{r,0}(\theta) \bigr\rrvert
^{2} \bigr]=\lim_{n}E_{0}
\frac{1}{n} \bigl\llvert M_{r,n}(\theta) \bigr\rrvert
^{2},
\end{equation}
$\mathbb{P}$-a.s. And (\ref{sigsquasyvar}) follows from (\ref{sigthe}) and
(\ref{expdrasanave}) because, under (\ref{appiner})
\begin{eqnarray*}
0&\leq&\limsup_{n} \biggl\llvert \biggl(E_{0}
\frac{1}{n} \bigl\llvert S_{n}(\theta)-E_{0}S_{n}(
\theta) \bigr\rrvert^{2} \biggr)^{1/2}-\sigma(\theta) \biggr
\rrvert
\\
&\leq&\limsup_{r}\limsup_{n} \biggl
\llvert \biggl(E_{0}\frac{1}{n} \bigl\llvert S_{n}(
\theta)-E_{0}S_{n}(\theta) \bigr\rrvert^{2}
\biggr)^{1/2}- \bigl(E \bigl\llvert D_{r,0}(\theta) \bigr
\rrvert ^{2} \bigr)^{1/2} \biggr\rrvert
\\
&\leq&\limsup_{r}\limsup_{n} \biggl(
\biggl\llvert \biggl(E_{0}\frac{1}{n} \bigl\llvert
S_{n}(\theta)-E_{0}S_{n}(\theta) \bigr\rrvert
^{2} \biggr)^{1/2}- \biggl(E_{0}\frac{1}{n}
\bigl\llvert M_{r,n}(\theta) \bigr\rrvert^{2}
\biggr)^{1/2} \biggr\rrvert\\
&&{}+ \biggl\llvert \biggl(E_{0}
\frac{1}{n} \bigl\llvert M_{r,n}(\theta) \bigr\rrvert
^{2} \biggr)^{1/2}- \bigl(E \bigl\llvert D_{r,0}(
\theta ) \bigr\rrvert^{2} \bigr)^{1/2} \biggr\rrvert \biggr)
\\
&=& \limsup_{r}\limsup_{n}\biggl( \biggl\llvert
\biggl(E_{0}\frac{1}{n} \bigl\llvert S_{n}(\theta
)-E_{0}S_{n}(\theta) \bigr\rrvert^{2}
\biggr)^{1/2}- \biggl(E_{0}\frac{1}{n} \bigl\llvert
M_{r,n}(\theta) \bigr\rrvert^{2} \biggr)^{1/2}
\biggr\rrvert\biggr)
\\
&\leq&\limsup_{r}\limsup_{n}
\biggl(E_{0}\frac{1}{n} \bigl\llvert S_{n}(\theta
)-E_{0}S_{n}(\theta)-M_{r,n}(\theta) \bigr
\rrvert^{2} \biggr)^{1/2}=0,
\end{eqnarray*}
$\mathbb{P}$-a.s.
\end{remark}
{\bf Proof of Lemma \ref{gengenbarinvpri}:} For $m\geq1$, the Skorohod
metric $d_{m}$ on $D[[0,m],\mathbb{C}]$ is dominated by the uniform
(product) metric. So
\[
d_{m}\bigl(r_{m}W_{n}(\theta),r_{m}V_{n}(
\theta)\bigr)\leq\frac{\sqrt{m}}{
\sqrt{n'}}\max_{1\leq k\leq n' }\bigl\llvert
S_{k}(\theta )-E_{0}S_{k}(\theta)-M
_{r,k}(\theta)\bigr\rrvert,
\]
where $n':=mn$ and $r_{m}$ is the restriction operator $D[[0,
\infty),\mathbb{C}]\to D[[0,m],\mathbb{C}]$. It follows from
(\ref{appiner}) that
\[
\lim_{r}\limsup_{n}\bigl\llVert
d_{m}\bigl(r_{m}W_{n}(\theta
),r_{m}V_{r,n}(\theta)\bigr)\bigr\rrVert_{
\mathbb{P}_{\omega},2}=0
\]
for $\mathbb{P}$-almost every $\omega$. Thus by an application of
Corollary~\ref{corlimlimlem} in the {appendix}, there exist random
functions $B^{m}_{r}(\theta)$ and $B^{m}_{\infty}(\theta)$ in
$D[[0,m],\mathbb{C}]$ such that the conclusion of Lemma~\ref{gengenbarinvpri} is satisfied with $B^{m}_{r}(\theta)$ in place
of $B_{r}(\theta)$ and with $B^{m}_{\infty}(\theta)$ in place of
$B(\theta)$. The first conclusion follows from a suitable adaptation
of \cite{Bilconpromea}, p.~173, Lemma~3 (or likewise adapting the
comments in 2. on Section~\ref{thespaD}).

To prove the second assertion, we argue as follows: first, for such
$\theta$'s, $B_{r}(\theta)$ is given by (\ref{limindbror}), and we
have seen that the (quenched) limit of (\ref{funqueinvpri}) as
$n\to\infty$ is the same as that of $B_{r}(\theta)$ as $r\to
\infty$.

Now, by convergence of Types theorems (see, for instance,
\cite{BilProMea}, p.~193 or \cite{baranexa}, Proposition~4),
$B_{r}(\theta)$ admits a limit in distribution if and only if there
exists
\[
\sigma^{2}(\theta):=\lim_{r}{E\bigl\llvert
D_{r,0}(\theta)\bigr\rrvert^{2}},
\]
which leads us to conclude that (\ref{funqueinvpri}) converges in the
quenched sense to
\[
\frac{\sigma(\theta)}{\sqrt{2}}\bigl(\mbox{\upshape Re}(B)+i\,\mbox{\upshape Im}(B)\bigr)
\]
as claimed.
\qed

\medskip

\begin{remark}
\label{conque}
Note that the convergence of (\ref{appfunr}) in the quenched sense as
$n\to\infty$ is guaranteed for all $\theta\in I$ and all
$0\leq r\leq\infty$, where $I$ is the set specified by Theorem~\ref{BaPeTh}. Whether the hypothesis (\ref{appiner}) is true as well for $\lambda$-a.e.
such $\theta$ is a question yet to be answered (it is  without the ``$\max_{1\leq k\leq n}$'', as can be seen from \cite{barpel}, p.~289, 290). Theorem~\ref{barinvprihancon} below gives particular cases of this statement.
\end{remark}

It is important also to point out that $(X_{k}-E_{-\infty}X_{k})$ is
a \textit{regular process}, according to the following
definition:

\medskip

\begin{dfn}[Regularity]
\label{regcondef}
Let $(Y_{k})_{k\in\mathbb{Z}}=(T^{k}Y_{0})_{k\in\mathbb{Z}}$ be a
$(\mathcal{F}_{k})_{k\in\mathbb{Z}}$-adapted stationary $p$-integrable
process ($p\geq1$). The process $(Y_{k})_{k\in\mathbb{Z}}$ is called
\textit{regular} if $E[Y_{0}|\mathcal{F}_{-k}]$ converges
to $0$ in $L^{p}_{\mathbb{P}}$. This is, if
%
\begin{equation}
\label{reg} \lim_{k} E\llvert E_{-k}Y_{0}
\rrvert^{p}=\lim_{k} E\llvert E_{0}Y_{k}
\rrvert^{p}=0.
\end{equation}
Equivalently, $(Y_{k})_{k\in\mathbb{Z}}$ is regular if  for
every $k\in\mathbb{Z}$
%
\begin{equation}
\label{regconequ} E[Y_{k}|\mathcal{F}_{-\infty}]=0.
\end{equation}
\end{dfn}

\begin{remark}
\label{remunireg}
As with the case of quenched convergence, the notion of regularity
depends on the choice of the $T$-filtration $(\mathcal{F}_{k})_{k
\in\mathbb{Z}}$, but it is easy to show that if $(Y_{k})_{k\in
\mathbb{Z}}$ is regular with respect to \textit{some} $T$-filtration
$(\mathcal{F}_{k}')_{k\in\mathbb{Z}}$ with $\sigma(Y_{0})\subset
\mathcal{F}_{0}'$, then it is regular with respect to its \textit{minimal
$T$-filtration} $(\mathcal{M}_{k})_{k\in\mathbb{Z}}=(T^{-k}
\mathcal{M}_{0})_{k\in\mathbb{Z}}$, where
\[
\mathcal{M}_{0}=\bigcap_{\alpha} \mathcal{G}_{\alpha},
\]
with the intersection running over the sigma algebras $\mathcal{G}_{
\alpha}$ in $\Omega$ for which $\sigma(Y)\subset\mathcal{G}_{
\alpha}\subset T^{-1}\mathcal{G}_{\alpha}$. In this sense, the notion
of regularity can be made independent of the choice of $\mathcal{F}
_{0}$.
\end{remark}

To prove the equivalence stated in Definition~\ref{regcondef} note that,
by (\ref{intkooopeequint}), (\ref{intkooopeequinf}), (\ref{revmarcon})
and the $L^{p}_{\mathbb{P}}$-continuity of the Koopman operator $T$, the
following equalities hold both $\mathbb{P}$-a.s. and in $L^{p}_{
\mathbb{P}}$ (with $p\geq1$): if $Y\in L^{p}_{\mathbb{P}}$
%
\begin{equation}
\label{equinvprotai} T^{k}E_{-\infty}Y=T^{k}\lim
_{j}E_{-j+k}Y=\lim_{j}T^{k}E_{-j+k}Y=
\lim_{j}E_{-j}T^{k}Y=E_{-\infty}T^{k}Y.
\end{equation}

Let us give now a result that, in the context of the present discussion,
will allow us to give a further characterization of the regularity
condition (\ref{reg}) (see Corollary~\ref{corcarreg}) and is of interest
by itself. We will use also this result later to prove part of the
statement of Theorem~\ref{Barinvpri} below.

\medskip

\begin{prop}
\label{proconsqu}
The equality
%
\begin{equation}
\label{seqergpow} \lim_{n}\frac{1}{n}\sum
_{k=0}^{n-1}\llvert E_{0}X_{k}
\rrvert^{2}=\lim_{n}\llVert E_{0}X
_{n}\rrVert_{\mathbb{P},2}^{2}
\end{equation}
holds $\mathbb{P}$-a.s.
\end{prop}

{\bf Proof:} $E_{0}T$ is a positive Dunford--Scwhartz operator and
therefore it is pointwise and mean ergodic (\cite{eisfarhaanag},
Theorem~8.24 and pp.~217--218): for every $Y\in L^{p}_{\mathbb{P}}$ ($p
\geq1$),
%
\begin{equation}
\label{ergthee0t} \frac{1}{n}\sum_{k=0}^{n-1}E_{0}T^{k}Y=E_{0}E[Y|\mathcal{T}] =EY,
\end{equation}
$\mathbb{P}$-a.s. and in $L^{p}_{\mathbb{P}}$, where $\mathcal{T}$ is
the invariant sigma algebra of $T$ and the second equality follows from
the assumption that $T$ is ergodic (we can identify the limit function
using the continuity of the operator $E_{0}$ in $L^{p}_{\mathbb{P}}$ and the mean
ergodic theorem for Koopman operators).

Now notice that, since for every $Z\in L^{2}_{\mathbb{P}}$ and
$(n,j)\in\mathbb{N}\times\mathbb{N}$,
$|E_{0}T^{n+j}Z|^{2}=|(E_{0}T)^{n+j}Z|^{2}
\leq E_{0}T^{n}|E_{0}T^{j}Z|^{2}$, $\mathbb{P}$-a.s., the process
$(F_{n})_{n\in\mathbb{N}}$ specified by
%
\begin{equation}
\label{deffn} F_{n}:=\sum_{k=0}^{n-1}E_{0}T^{k}
\llvert X_{0}\rrvert^{2}-\sum_{k=0}^{n-1}
\bigl\llvert E_{0}T ^{k}X_{0} \bigr\rrvert
^{2}
\end{equation}
is nonnegative and superadditive (\cite{Kre}, p.~146., or
\cite{akcsuc}).

An application of \cite{Kre}, Theorem~5.3, together with the
Dunford-Schwartz pointwise ergodic theorem, implies that there exists
a (nonnegative) function $Y\in L^{1}_{\mathbb{P}}$ such that
\begin{enumerate}
\item[1.] $EY=\lim_{n}EF_{n}/n$,
\item[2.] For every $n\in\mathbb{N}$
\[
F_{n}\leq G_{n}:=\sum_{k=0}^{n-1}E_{0}T^{k}Y,
\]
$\mathbb{P}$-a.s., and
\item[3.] $\lim_{n}{F_{n}}/{G_{n}}=1$, $\mathbb{P}$-a.s. on $
\bigcup_{k\in\mathbb{N}}[E_{0}T^{k}Y>0]$ (either this set or its
complement has probability one in this case, but this is not assumed for
what follows along this proof).
\end{enumerate}
It follows from these observations (note also that $F_{n}=0$,
$\mathbb{P}$-a.s., on $\bigcap_{k\in\mathbb{N}}[E_{0}T^{k}Y=0]$) and
another application of the Dunford-Schwartz ergodic theorem that
%
\begin{equation}
\label{equprolasrep} \lim_{n}\frac{1}{n}\sum
_{k=0}^{n-1} \bigl\llvert E_{0}T^{k}X_{0}
\bigr\rrvert^{2}=E \bigl[\llvert X_{0}\rrvert
^{2}-Y \bigr],
\end{equation}
$\mathbb{P}$-a.s.\footnote{An alternative argument is the following:
by (\ref{ergthee0t}), the invariant sigma algebra of the operator
$E_{0}T$ is trivial, and it follows by an application of Corollary~1 in
\cite{akcsuc} with $s_{k}'=k$ that ${F_{n}}/{n}\to_{n} EY$,
$\mathbb{P}$-a.s. The pointwise ergodicity of $E_{0}T$ (i.e., a novel
application of (\ref{ergthee0t})) gives (\ref{equprolasrep}).}

To actually show that $EY=\|X_{0}\|_{\mathbb{P},2}^{2}-\lim_{n}\|E
_{0}T^{n}X_{0}\|_{\mathbb{P},2}^{2}$ we first notice that
\[
EY=\lim_{n}\frac{EF_{n}}{n}=\lim_{n}
\Biggl(\llVert X_{0}\rrVert_{\mathbb{P},2}^{2}-
\frac{1}{n}\sum_{k=0}^{n-1}\bigl\llVert
E_{0}T^{k}X_{0}\bigr\rrVert_{\mathbb{P},2}^{2}
\Biggr)
\]
$\mathbb{P}$-a.s., and that by Jensen's inequality
\[\bigl\|E_{0}T^{k+1}X_{0}\bigr\|_{
\mathbb{P},2}^{2}:=E\bigl|E_{0}T^{k+1}X_{0}\bigr|^{2}\leq
E\bigl|E_{0}T^{k}X_{0}\bigr|^{2}=:\bigl\|E
_{0}T^{k}X_{0}\bigr\|_{\mathbb{P},2}^{2}.\]
Thus the sequence
$(\|E_{0}T^{k}X_{0}\|_{\mathbb{P},2}^{2})_{k\in\mathbb{N}}$ is
decreasing. It follows that
\[
\lim_{n}\frac{1}{n}\sum_{k=0}^{n-1}
\bigl\llVert E_{0}T^{k}X_{0}\bigr\rrVert
_{\mathbb{P},2} ^{2}=\lim_{n}\bigl\llVert
E_{0}T^{n}X_{0}\bigr\rrVert_{\mathbb{P},2}^{2}.
\]\qed
%

\begin{remark}
\label{remconpow}
Using exactly the same argument as in the proof of Proposition~\ref{proconsqu}, one can further prove the $\mathbb{P}$-a.s. convergence
\[
\lim_{n}\frac{1}{n}\sum_{k=0}^{n-1}
\bigl\llvert E_{0}T^{k}X_{0}\bigr\rrvert
^{\beta}=\lim_{n}E\llvert E_{0}X_{n}
\rrvert^{\beta}
\]
for every $1\leq\beta\leq p$, provided that $X_{0}\in L^{p}_{
\mathbb{P}}$. A~similar argument proves the same conclusion when
$0\leq\beta< 1\leq p$ (in this case the process $(\sum_{k=0}^{n-1}|E
_{0}T^{k}X_{0}|^{\beta})_{n\in\mathbb{N}^{*}}$ is
superadditive).\footnote{There is an additional detail in this case:
to prove that the increasing, Ces\`{a}ro-convergent sequence
$(E|E_{0}T^{n}X|^{\beta})_{n\in\mathbb{N}}$ is bounded (i.e.,
convergent) note for instance that, by H\"{o}lder's inequality and
Jensen's inequality
\[
\bigl\llVert E_{0}T^{n}X\bigr\rrVert_{\mathbb{P},\beta}\leq
\bigl\llVert E_{0}T^{n}X\bigr\rrVert_{
\mathbb{P},p}\leq
\llVert X\rrVert_{\mathbb{P},p}
\]
for all $n\in\mathbb{N}$ (it is actually a short exercise to prove that
an unbounded and increasing sequence is \textit{not} Ces\`{a}ro convergent:
its Ces\`{a}ro averages also diverge to $\infty$).}

More generally, an easy adaptation of the first part of the arguments
in this proof (using the same results from \cite{Kre}) proves that if
$L$ is a positive Dunford-Schwartz operator (or any conservative
positive contraction that is also pointwise ergodic) and $\varphi:[0,
\infty)\to[0,\infty)$ is a nonnegative measurable function for which
$\varphi\circ LX\leq L(\varphi\circ X)$, $\mathbb{P}$-a.s. for a given
nonnegative $X\in L^{p}_{\mathbb{P}}$ ($p\geq1$), then
\[
\frac{1}{n}\sum_{k=0}^{n-1}\varphi
\bigl(L^{k}X\bigr)
\]
converges $\mathbb{P}$-a.s.
\end{remark}

The following is an obvious consequence of Proposition~\ref{proconsqu}.

\medskip

\begin{cor}
\label{corcarreg}
The process $(X_{k})_{k\in\mathbb{Z}}$ is regular (Definition~\ref{regcondef}) if and only if
\[
\lim_{n}\frac{1}{n}\sum_{k=0}^{n-1}
\llvert E_{0}X_{k}\rrvert^{2}=0,
\]
$\mathbb{P}$-a.s.
\end{cor}

Let us now introduce is the ``averaged-frequency'' version of Lemma~\ref{convnr}, which will be the building block for the proof of Theorem~\ref{Barinvpri}.

\medskip
\begin{lemma}
\label{convnrdou}
If $\mathcal{B}_{0}=\{\emptyset,[0,2\pi)\}$ is the trivial sigma
algebra and $0\leq r\leq\infty$ is fixed, the random function
$(\theta,\omega)\mapsto V_{r,n}(\theta,\omega)$ defined by (\ref{appfunr}) (where $V_{\infty
,n}(\theta,\omega):=0$ if $\theta\notin I$) converges in the quenched
sense with respect to $\mathcal{B}_{0}\otimes\mathcal{F}_{0}$ (see
Section~\ref{queconprospa}), as $n\to\infty$, to the function
$(\theta,\omega')\mapsto B_{r}(\theta,\omega')$ defined by (\ref{limindbror}) (where $B(\theta
,\omega')=0$ if $\theta\notin I$). Equivalently, there exists
$\Omega_{0}\subset\Omega$ with $\mathbb{P}\Omega_{0}=1$ such that for
every $\omega\in\Omega_{0}$, $V_{r,n}\Rightarrow_{n} B_{r}$ under
$\lambda\times\mathbb{P}_{\omega}$.
\end{lemma}

{\bf Proof:} This follows at once from Corollary~\ref{marcasdou} in
Section~\ref{secmarcas}.\qed

We will use Lemma~\ref{convnrdou} in order to prove Theorem~\ref{Barinvpri} below. We will also need the following lemma, which is
essentially the same as Lemma~\ref{gengenbarinvpri} except that, just
as in Lemma~\ref{convnrdou}, the processes involved are seen as random
elements whose domain is the product space $[0,2\pi)\times\Omega$,
and the notion of quenched convergence is understood with respect to
$\mathcal{B}_{0}\otimes\mathcal{F}_{0}$. The proof is exactly as that
of Lemma~\ref{gengenbarinvpri}.

\medskip

\begin{lemma}
\label{gengenbarinvpridou}
With the notation introduced in Section~\ref{queconprospa} and in (\ref{firmardifr})-(\ref{funqueinvprinoncen}) (see also
Remark~\ref{remdouint}), and taking $\mathcal{B}_{0}:=\{\emptyset,[0,2
\pi)\}$, assume that
%
\begin{equation}
\label{appinerdou} \lim_{r}\limsup_{n}
\tilde{E}_{0} \biggl[\frac{1}{n}\max_{1\leq k\leq n}
\bigl\llvert S _{k}(\theta,\omega)-E_{0}S_{k}(
\theta,\omega)-M_{r,k}(\theta, \omega) \bigr\rrvert^{2}
\biggr]=0
\end{equation}
$\mathbb{P}$-a.s. Then the random function $(\theta,\omega)\mapsto W_{n}(\theta,
\omega)$ specified by (\ref{funqueinvpri}) converges in the quenched sense as $n\to\infty$, with respect to $
\mathcal{B}_{0}\otimes\mathcal{F}_{0}$, to the random function $(\theta,\omega')\mapsto
B(\theta,\omega')$ specified by (\ref{limindbro}) ($B(\theta,\omega)=0$ if $\theta\notin I$). Equivalently, there exists
$\Omega_{0}\subset\Omega$ with $\mathbb{P}\Omega_{0}=1$ such that for
every $\omega\in\Omega_{0}$, $W_{n}\Rightarrow_{n} B$ under $\lambda
\times\mathbb{P}_{\omega}$.
\end{lemma}

\subsection{Main results}

We will prove that the hypotheses in Lemma~\ref{gengenbarinvpridou} hold
under the general setting established in Section~\ref{genset}. Let us
state the respective result in terms of the regular conditional measures
associated to $E_{0}$:


\medskip

\begin{thm}
\label{Barinvpri}
There exists $\Omega_{0}\subset\Omega$ with $\mathbb{P}\Omega_{0}=1$
such that for every $\omega_{0}\in\Omega_{0}$ the random functions
\[
W_{n}:\bigl([0,2\pi)\times\Omega, \mathcal{B}\otimes{
\mathcal{F}}, \lambda\times\mathbb{P}_{\omega_{0}}\bigr)\to D\bigl[[0,\infty),\mathbb{C}\bigr]
\]
specified by (\ref{funqueinvpri}) (see also Remark~\ref{remtwopar})
converge in distribution, as $n\to\infty$, to the random function
\[
B:\bigl([0,2\pi)\times\Omega',\mathcal{B}\otimes
\mathcal{F}',\lambda \times\mathbb{P}'\bigr)\to D\bigl[[0,\infty),
\mathbb{C}\bigr]
\]
specified by (\ref{limindbro}) (with $B(\theta,\omega)=0$ if $\theta\notin
I$). If (\ref{reg}) holds, then the non-centered cadlag functions
(\ref{funqueinvprinoncen}) satisfy the same conclusion.
\end{thm}

\medskip

Note that, by the argument preceding Remark~\ref{queimpsemque}, this
theorem implies Proposition~2.1 in \cite{PeWu}. This theorem will be
proved in Section~\ref{prothe}.

We will now present the cases in which the estimate (\ref{appiner}) is
proved along this paper. The notation introduced in the lemmas above is
kept here.

To present our results start by considering the following conditions of
weak dependence
\begin{enumerate}
\item[(a)] The \textit{generalized Hannan condition}, given by
%
\begin{equation}
\label{imhancon} \sum_{n\geq0} \bigl\llVert
\mathcal{P}_{0}(X_{n+1}-X_{n}) \bigr\rrVert
_{\mathbb{P},2}< \infty,
\end{equation}
which is clearly a ``weak'' version of the \textit{Hannan}
condition\footnote{To see that (\ref{imhancon}) is strictly weaker
than (\ref{hancon}) consider the process
\[
X_{k}:=\sum_{j\in\mathbb{N}^{*}}\frac{1}{j}x_{k-j},
\]
where $(x_{j})_{j\in\mathbb{Z}}$ are the coordinate functions in
$\mathbb{R}^{\mathbb{Z}}$, seen as an i.i.d. centered sequence in
$L^{2}$, $T$ is the left shift, and $\mathcal{F}_{0}=\sigma(x_{k})_{k
\leq0}$.}
%
\begin{equation}
\label{hancon} \sum_{n}\llVert
\mathcal{P}_{0}X_{n}\rrVert_{\mathbb{P},2}<\infty.
\end{equation}
Note, on the other side, that (\ref{imhancon}) is the same as
(\ref{hancon}) for the (stationary) process $(Y_{n})_{n}$ specified by
\[
Y_{n}:=X_{n}-X_{n-1}.
\]
\item[(b)] The \textit{Maxwell and Woodroofe condition for a fixed
frequency.} Which states that, for $\theta\in[0,2\pi)$,
%
\begin{equation}
\label{maxwoo} \sum_{k=1}^{\infty}
\frac{\|E_{0}S_{k}(\theta)\|_{\mathbb{P},2}}{k
^{3/2}}<\infty.
\end{equation}
Condition (\ref{maxwoo}) is a ``rotated'' version of the \textit{Maxwell
and Woodroofe condition}
%
\begin{equation}
\label{maxwoononrot} \sum_{k=1}^{\infty}
\frac{\|E_{0}S_{k}\|_{\mathbb{P},2}}{k^{3/2}}< \infty.
\end{equation}
\end{enumerate}

Conditions (\ref{imhancon}) and (\ref{maxwoononrot}) are logically
independent, see \cite{duri},\hskip.2pt\footnote{More precisely, the results
stated in \cite{duri} give that (\ref{imhancon}) is not sufficient for
(\ref{maxwoononrot}). To prove that (\ref{imhancon}) is not necessary
for (\ref{maxwoononrot}) consider the process constructed in
\cite{duri}, and strengthen the assumption ``$N_{k+1}>N_{k}$'' to ``$N
_{k+1}>N_{k}+1$''. Following the arguments in that paper, it is easy to
see that (\ref{imhancon}) and (\ref{hancon}) are equivalent for the
process under consideration, and the proof of independence can be
carried over by substituting, in the first counterexample on
\cite{duri}, ``$N_{k}=k$'' by ``$N_{k}=3k$''.} and they imply the
existence of martingale approximations allowing us to prove the
following theorem:

\medskip
 
\begin{thm}
\label{barinvprihancon}
With the notation introduced in lemmas \ref{convnr} and
\ref{gengenbarinvpri}, and assuming (\ref{imhancon}) or (\ref{maxwoo}),
if $e^{2i\theta}\notin \mathrm{Spec}_{p}(T)$, then (\ref{funqueinvpri})
converges in the quenched sense, as $n\to\infty$, to
%
\begin{equation}
\label{quelimfixfreequ} \omega'\mapsto\frac{\sigma(\theta)}{\sqrt
{2}} \bigl(
\mbox{\upshape Re}(B) \bigl(\omega' \bigr)+i\,\mbox{\upshape Im}(B) \bigl( \omega
' \bigr) \bigr).
\end{equation}
If (\ref{maxwoo}) holds, (also) the non-centered cadlag functions
(\ref{funqueinvprinoncen}) converge in the quenched sense, as
$n\to\infty$, to (\ref{quelimfixfreequ}).
\end{thm}

It is important to point out the following distinction implicit in the
statement of Theorem~\ref{barinvprihancon}: when condition
(\ref{imhancon}) holds, Theorem~\ref{barinvprihancon} tells us that,
\textit{for every} $\theta$ with $e^{2i\theta}\notin \mathrm{Spec}_{p}(T)$, and
in particular for $\lambda$-a.e. $\theta$ (see
(\ref{lammeapoispezer})), (\ref{funqueinvpri}) converges in the quenched
sense, whereas the conclusion under the Maxwell and Woodroofe condition
(\ref{maxwoo}) is that for \textit{a given} $\theta$ with $e^{2i\theta}\notin \mathrm{Spec}
_{p}(T)$, (\ref{funqueinvpri}) \textit{and} (\ref{funqueinvprinoncen})
converge in the quenched sense (to the same limit function) \textit{if}
(\ref{maxwoo}) holds.

It is therefore desirable to give a criterion implying the fulfillment
of (\ref{maxwoo}) for $\lambda$-a.e. $\theta$. This is the content of
the following result:

\medskip
\begin{thm}
\label{promaxwooae}
If for some $\beta>1$
%
\begin{equation}
\label{equpromaxwooae} \sum_{k=1}^{\infty} (\log
k)^{\beta}\frac{\|E_{0}X_{k}\|_{
\mathbb{P},2}^{2}}{k}<\infty,
\end{equation}
then (\ref{maxwoo}) holds for $\lambda$-a.e. $\theta$. In particular
(\ref{equpromaxwooae}) implies that (\ref{funqueinvprinoncen}) (and
(\ref{funqueinvpri})) converges in the quenched sense, as $n\to
\infty$, for $\lambda$-a.e. $\theta$, as specified in Theorem~\ref{barinvprihancon}.
\end{thm}

\medskip
\begin{remark}
\label{remrancen}
One may naturally ask if, just as in the case of Theorem~\ref{Barinvpri}, the verification of (\ref{reg}) \textit{together} with the
quenched convergence of (\ref{funqueinvpri}) for a given $\theta$ (or
for $\lambda$-a.e. $\theta$) are sufficient conditions for the
quenched convergence of (\ref{funqueinvprinoncen}) as $n\to\infty$.

The actual answer is \textit{no}: consider one more time $(x_{j})_{j
\in\mathbb{Z}}$, $(\mathcal{F}_{j})_{j\in\mathbb{Z}}$, $T$ as in the
footnote following (\ref{imhancon}) and note again that any (centered)
linear process
%
\begin{equation}
\label{linpro} X_{k}=\sum_{j\in\mathbb{N}}a_{j}x_{k-j}
\end{equation}
with $(a_{j})_{j\in\mathbb{N}}\in l^{1}({\mathbb{N}})$ ($(x_{j})_{j
\in\mathbb{N}}$ i.i.d., $Ex_{0}=0$, $E|x_{0}|^{2}=1$) satisfies
(\ref{hancon}) and therefore also (\ref{imhancon}). Thus the conclusion
of Theorem~\ref{barinvprihancon} holds for (\ref{funqueinvpri}) when
$(X_{k})_{k\in\mathbb{Z}}$ is given by (\ref{linpro}), $(a_{j})_{j
\in\mathbb{N}}\in l^{1}({\mathbb{N}})$ and $\theta\in(0,2 \pi)
\setminus\{\pi\}$ (see Corollary~\ref{weamixcas} below).

Such a process is necessarily regular because, by Kolmogorov's zero-one
law, $E_{-\infty}Z=EZ$, $\mathbb{P}$-a.s. for every $Z\in L^{1}_{
\mathbb{P}}$ in this case but, as shown in \cite{baranexa}, the coefficients
$(a_{k})_{k\in\mathbb{N}}$ can be chosen in such a way that, for every
$\theta\in[0,2\pi)$, $S_{n}(\theta)/\sqrt{n}$ admits no limit in
distribution under $\mathbb{P}_{\omega}$ for $\mathbb{P}$-a.e.
$\omega$. This is of course an obstruction (consider $t=1$) to the
quenched convergence, as $n\to\infty$, of (\ref{funqueinvprinoncen}) to
(\ref{quelimfixfreequ})  or to any random function that is
continuous at $1$ with probability one, but we can discard convergence
to any other random function in $D[[0,\infty),\mathbb{C}]$ considering
the identity $U_{n}(\theta,t)=( \lfloor nt \rfloor/n)^{1/2}U
_{ \lfloor nt \rfloor}(\theta,1)$ ($U_{n}$ is given by
(\ref{funqueinvprinoncen})) and the general arguments behind the
observation 2. in Section~\ref{thespaD}. See for instance
\cite{Bilconpromea}, p.~138.

These considerations support the claim that (even) under the assumptions
presented along this paper, including the regularity condition
(\ref{reg}), the quenched functional convergence with respect to
$\mathcal{F}_{0}$ of $S_{n}(\theta)/\sqrt{n}$ for $\lambda$-a.e.
{fixed} $\theta$ as $n\to\infty$ is a conclusion \textit{strictly
stronger} than that of the quenched functional convergence with respect
to $\{\emptyset,[0,2\pi)\}\otimes\mathcal{F}_{0}$ of the same process
seen as a process parametrized by $(\theta,\omega)$. More explicitly,
if $(X_{k})_{k\in\mathbb{Z}}$ is the process constructed in
\cite{baranexa} then, by Theorem~\ref{Barinvpri}, $(\theta,\omega)
\mapsto U_{n}(\theta,\omega)$ converges in the quenched sense with
respect to $\{\emptyset,[0,2\pi)\}\otimes\mathcal{F}_{0}$ to
$(\theta,\omega')\mapsto B(\theta,\omega')$, and we have just
mentioned that (\ref{funqueinvprinoncen}) does not converge (in
$D[[0,\infty),\mathbb{C}]$) in the quenched sense with respect to
$\mathcal{F}_{0}$ for any fixed $\theta\in[0,2\pi)$, so that we
cannot obtain the conclusion of Theorem~\ref{barinvprihancon} for
(\ref{funqueinvprinoncen}) from the regularity condition (\ref{reg})
\textit{in spite that} the (full) conclusion of Theorem~\ref{Barinvpri} is true
for a process satisfying this condition. {The question of whether
this is also the case for \textup{(\ref{funqueinvpri})} seems to be still open},
as hinted already in Remark~\ref{conque}.

Finally, note that these kind of considerations show that one cannot
deduce the fulfillment of (\ref{maxwoo}) for $\lambda$-a.e.
$\theta$, not even for any given fixed $\theta$, from the fulfillment
of (\ref{imhancon}) only.\footnote{Except perhaps, though one would
not expect such thing, if $e^{2i\theta}\in \mathrm{Spec}_{p}(T)$ and
$\theta\neq0$ are conditions imposed on $\theta$. The case
$\theta=0$ being covered by \cite{duri} as specified above, and the
case $e^{2i\theta}\notin \mathrm{Spec}_{p}(T)$ being covered by our
considerations.} In particular, (\ref{equpromaxwooae}) does not follow
from (\ref{imhancon}).
\end{remark}

\medskip

\begin{remark}
\label{remhanaereg}
It is also convenient to point out the relation between
(\ref{equpromaxwooae}) and (\ref{reg}): first, it clearly follows from
(\ref{equpromaxwooae}) that\vspace*{-2pt}
%
\begin{equation}
\label{almreg} \sum_{k=1}^{\infty}
\frac{\|E_{0}X_{k}\|_{\mathbb{P},2}^{2}}{k}< \infty
\end{equation}
and this, together with Kronecker's lemma, implies that under
(\ref{equpromaxwooae})\vspace*{-2pt}
\[
\lim_{n}\frac{1}{n}\sum_{k=0}^{n-1}
\llVert E_{0}X_{k}\rrVert_{\mathbb{P},2} ^{2}=0,
\]
which is the same as (\ref{reg}) because $(\|E_{0}X_{k}\|_{
\mathbb{P},2}^{2})_{k\in\mathbb{N}}$ is decreasing.\vspace*{1pt}

On the other side, (\ref{reg})  does not imply (\ref{equpromaxwooae})
because, as we have seen in Remark~\ref{remrancen}, there exist
processes satisfying (\ref{reg}) that do not satisfy the conclusion of
Theorem~\ref{barinvprihancon}. Thus, (\ref{equpromaxwooae}) is \textit
{strictly stronger} than (\ref{reg}).

Collecting our results, the following question deserves to be pointed
out: \textit{is there $0\leq\beta\leq1$ such that
\textup{(\ref{equpromaxwooae})} holds for a process that does} not
\textit{satisfy the conclusion of Theorem \textup{\ref{barinvprihancon}}? (with and
without the statement of convergence for \textup{(\ref{funqueinvprinoncen})})}.
And if the answer is ``no'', \textit{what conditions ``in between''
\textup{(\ref{almreg})} and \textup{(\ref{reg})} break down the conclusion of Theorem
\textup{\ref{barinvprihancon}}?}
\end{remark}

To finish this section, it is worth to further specify a case in which
the set of frequencies where the asymptotic distribution is as in
(\ref{quelimfixfreequ}) can be easily described. To motivate the
following result, recall that $T$ is weakly mixing if and only if
$\mathrm{Spec}_{p}(T)=\{1\}$ (see \cite{qua}, Section~8 for a review of this
and other related facts).

Now, as a subgroup of $\mathbb{T}$, $\mathrm{Spec}_{p}(T)$ is finite (actually:
closed) if and only if there exists $m\in\mathbb{N}^{*}$ such that\vspace*{-2pt}
%
\begin{equation}
\label{eigangfinspe} \mathrm{Spec}_{p}(T):= \bigl\{e^{2\pi k i/m} \bigr
\}_{k=0}^{m-1}.
\end{equation}
In other words, $\mathrm{Spec}_{p}(T)$ is finite if and only it it consists of
the points in the unit circle given by the {rational rotations} by an
angle of $2\pi/m$ or, what is the same, by the $m$th roots of unity.

\medskip
\begin{cor}
\label{weamixcas}
Assume that $\mathrm{Spec}_{p}(T)$ is finite and its elements are the $m$th
roots of unity. Under the hypotheses in Theorem~\ref{barinvprihancon},
(\ref{funqueinvpri}) converges in the quenched sense as $n\to\infty
$ to (\ref{quelimfixfreequ}) for all $\theta\in[0,2\pi)$ such that
$e^{2im\theta}\neq1$. If $T$ is in particular weakly mixing,
(\ref{quelimfixfreequ}) describes the asymptotic quenched limit of
(\ref{funqueinvpri}) for all $\theta\neq0,\pi$. The same conclusions
hold for (\ref{funqueinvprinoncen}) when (\ref{maxwoo}) holds.\looseness=-1\vadjust{\goodbreak}
\end{cor}

{\bf Proof:} Immediate from (\ref{eigangfinspe}), Theorem~\ref{barinvprihancon} and Theorem~\ref{promaxwooae}.\qed

\medskip
\begin{remark}
Of course, the conclusion of Theorem~\ref{barinvprihancon} is also true
for the weakly mixing case under (\ref{equpromaxwooae}) ($\beta>1$),
but our proof in this case losses track of the ($\lambda$-measure one)
set of frequencies $\theta$ where (\ref{maxwoo}) holds.
\end{remark}


\section{Martingale case}\label{secmarcas}
Since our proofs are based on martingale approximations we will study,
in this section, the asymptotic distributions of (\ref{appfunr}) for the
case in which $D_{r,0}(\theta)$ is replaced by an abstract martingale
difference $D_{0}(\theta)$.

The result is Theorem~\ref{marcas}. The martingale case on the product
space (``averaged frequency'' case) follows from this one via some of
the general results presented in this paper, and it is the content of
Corollary~\ref{marcasdou}.


\medskip

\begin{thm}
\label{marcas}
Under the setting in Section~\ref{secass}, and with the notation in
Lemma~\ref{convnr}, given $\theta\in[0,2\pi)$ such that $e^{2i
\theta}\notin \mathrm{Spec}_{p}(T)$, assume that $D_{0}(\theta)\in L^{2}_{
\mathbb{P}}(\mathcal{F}_{0})\ominus L^{2}_{\mathbb{P}}(\mathcal{F}
_{-1})$ is given, and define the $(\mathcal{F}_{k-1})_{k\in
\mathbb{N}^{*}}$-adapted martingale $(M_{k}(\theta))_{k\in
\mathbb{N}}$ by
%
\begin{equation}
\label{nmardef} M_{n}(\theta):=\sum_{k=0}^{n-1}T^{k}D_{0}(
\theta)e^{ik\theta}
\end{equation}
for all $n\in\mathbb{N}$. Then the sequence $(V_{k}(\theta))_{k
\in\mathbb{N}^{*}}$ of random elements of $D[[0,\infty),\mathbb{C}]$
defined by
%
\begin{equation}
\label{defofvnthe} V_{n}(\theta) (t):= M_{ \lfloor nt \rfloor
}(\theta)/
\sqrt{n}
\end{equation}
for every $n\in\mathbb{N}^{*}$, converges in the quenched sense with
respect to $\mathcal{F}_{0}$ to the random function $B(\theta):
\Omega'\to D[[0,\infty),\mathbb{C}]$ given by
%
\begin{equation}
\label{asyvnthemarcas} B(\theta) \bigl(\omega' \bigr)= \bigl[E \bigl\llvert
D_{0}(\theta) \bigr\rrvert^{2}/2 \bigr]^{1/2}
\bigl( \mbox{\upshape Re}(B) \bigl(\omega' \bigr)+i\,\mbox{\upshape Im}(B) \bigl( \omega
' \bigr) \bigr).
\end{equation}
\end{thm}

\medskip
\begin{remark}
\label{remequmarcas}
Before proceeding to the proof it is worth noticing the following: the
conclusion of Theorem~\ref{BaPeTh}, specialized to this case, is a
statement about the asymptotic distribution of the random variables
$V_{n}(\theta)(1)$. By Lemma~\ref{dedmerpel} in the {appendix} and the
orthogonality under $E_{0}$ of $(T^{k}D_{0}(\theta))_{k\in
\mathbb{N}}$,\footnote{Note that if $(k,r)\in\mathbb{N}\times
\mathbb{N}^{*}$ is given then, since $T^{r}{D_{0}(\theta)}\in L^{2}
_{\mathbb{P}}(\mathcal{F}_{r})\ominus L^{2}_{\mathbb{P}}(\mathcal{F}
_{r-1})$,
\[
E_{0}\bigl[T^{k}D_{0}(\theta)T^{k+r}
\overline{D_{0}(\theta)}\bigr]=E_{0}\bigl[T
^{k}D_{0}(\theta)E_{k}T^{k+r}
\overline{D_{0}(\theta)}\bigr]=E_{0}T^{k}\bigl[D
_{0}(\theta)E_{0}T^{r}\overline{D_{0}(
\theta)}\bigr]=0,
\]
$\mathbb{P}$-a.s.}
\[
E\bigl[\bigl\llvert D_{0}(\theta)\bigr\rrvert^{2}\bigr]=
\lim_{n}\frac{1}{n}\sum_{k=1}^{n-1}E_{0}T^{k}
\bigl\llvert D _{0}(\theta)\bigr\rrvert^{2}=\lim
_{n}\frac{1}{n}E_{0}\bigl\llvert
M_{n}(\theta )-E_{0}M_{n}( \theta)\bigr
\rrvert^{2},
\]
and therefore the equality (\ref{asyvar}) is certainly verified in this case.\vadjust{\goodbreak}
\end{remark}

{\bf Proof of Theorem \ref{marcas}:} Let us start by sketching the
argument of the proof: we will see that for $\theta$ as specified,
there exists $\Omega_{\theta}\subset\Omega$ with $\mathbb{P}
\Omega_{\theta}=1$ such that for every $\omega\in\Omega_{\theta}$
the following holds:

\begin{enumerate}
\item[(a)]\textit{The sequence of random functions
$(V_{n}(\theta))_{n}$ in $D[[0,\infty),\mathbb{C}]$ is tight with
respect to $\mathbb{P}_{\omega}$.} To prove this, we will actually
prove the convergence in distribution of both the real and imaginary
parts of $(V_{n}(\theta))_{n}$ to a Brownian motion via Theorem~\ref{bilinvpri}.
\item[(b)] \textit{The finite dimensional asymptotic distributions under
$\mathbb{P}_{\omega}$ of $(V_{n}(\theta))_{n}$ converge\vspace*{1pt} to those of
two independent standard motions with the rescaling $E[(D_{0}(\theta
))^{2}]^{1/2}/
\sqrt{2}$.}  For this, we will proceed via
the Cram\'{e}r--Wold theorem, using some of the results already
presented.
\end{enumerate}

We go now to the details: first, we will assume, making it explicit only
when necessary, that $E_{0}$ is the version of $E[\,\cdot\,|
\mathcal{F}_{0}]$ given by integration with respect to the decomposing
probability measures $\{\mathbb{P}_{\omega}\}_{\omega\in\Omega}$.

Now denote, for every $k\in\mathbb{N}$
%
\begin{equation}
\label{defdkmarcas} D_{k}(\theta):=T^{k}D_{0}(\theta
).
\end{equation}

Let $\Omega_{\theta,1}'$ be the set of probability one guaranteed by
Lemma~\ref{rel16cunmerpel} for the case $Y=D_{0}(\theta)$. By Remark~\ref{remproone}, there exists a set $\Omega_{\theta,1}$ with
$\mathbb{P}\Omega_{\theta,1}=1$ such that for every $\omega\in
\Omega_{\theta,1}$
\[
\lim_{n}\frac{1}{n}\sum_{k=0}^{n-1}
E_{k-1}\bigl(z\cdot\bigl(D_{k}(\theta)e ^{ik\theta}\bigr)
\bigr)^{2}= \frac{\llvert z\rrvert^{2}}{2}E\bigl\llvert D_{0}(\theta)
\bigr\rrvert^{2}
\]
$\mathbb{P}_{\omega}$-a.s. for all $z\in\mathbb{C}$.

For such $\omega$'s the first hypothesis of Theorem \ref{bilinvpri} is verified  by the triangular arrays $(Re(M_{k}(\theta)/\sqrt{n}))_{1\leq k \leq n}$  and $(Im(M_{k}(\theta)/\sqrt{n}))_{1\leq k \leq n}$ ($n\in\mathbb{N}^{*}$) with respect to $\mathbb{P}_{\omega}$, because they arise from the particular choices $z=1$ and $z=i$ respectively.

To verify the second hypothesis in Theorem \ref{bilinvpri} we start from the $\mathbb{P}-$a.s. inequality

To verify the second hypothesis in Theorem~\ref{bilinvpri}, we start
from the $\mathbb{P}$-a.s. inequality
%

$$E_{0}\left[\frac{1}{n}\sum_{k=0}^{n-1} \left((Re(D_{k}(\theta)e^{ik\theta}))^2I_{[|Re(D_{k}(\theta)e^{ik\theta})|\geq \epsilon \sqrt{n}]]}+(Im(D_{k}(\theta)e^{ik\theta}))^2I_{[|Im(D_{k}(\theta)e^{ik\theta})|\geq \epsilon\sqrt{n}]}\right)\right] $$
\begin{equation}
\label{inetigreaima}
\leq E_{0}\left[\frac{1}{n}\sum_{k=0}^{n-1} |D_{k}(\theta)|^2I_{[|D_{k}(\theta)|\geq \epsilon\sqrt{n}]}\right].
\end{equation}

Now, given $\eta>0$ there exists $N\geq 0$ such that $\mu_{N}:= {E}[|D_{0}(\theta)|^{2}I_{[|D_{0}(\theta)|^{2}\geq \epsilon^{2} N]}]<\eta$, and therefore
$$\limsup_{n}\frac{1}{n}\sum_{k= 0}^{n-1}E_{0}T^{k}[|D_{0}(\theta)|^{2}I_{[|D_{0}(\theta)|^{2}\geq \epsilon^{2} n]}]\leq$$
\begin{equation}
\label{inelimequzer}
\limsup_{n}\frac{1}{n}\sum_{k=0}^{ n-1}E_{0}T^{k}[|D_{0}(\theta)|^{2}I_{[|D_{0}(\theta)|^{2}\geq \epsilon^{2} N]}]=\mu_{N}\leq \eta
\end{equation}
over a set $\Omega_{\theta,\epsilon,\eta}$ with $\mathbb{P}
\Omega_{\theta,\epsilon,\eta}=1$, where we made use of Corollary~\ref{dedmerpel}. Without loss of generality, (\ref{inetigreaima}) holds
for all $\omega\in\Omega_{\theta,\epsilon,\eta}$.

Denote by $Z_{n}^{\epsilon}$ the random variable at the left-hand side
of the inequality (\ref{inetigreaima}) and note that, if we define
%
\begin{equation}
\label{defome2pri} \Omega_{\theta,2}=\bigcap_{\epsilon>0,\eta>0}
\Omega_{\theta,\epsilon,\eta}
\end{equation}
where the intersection runs over rational $\epsilon$, $\eta$, then
$\mathbb{P}\Omega_{\theta,2}=1$, and for every $\epsilon>0$ and every
$\omega\in\Omega_{\theta,2}$
\[
\lim_{n}Z_{n}^{\epsilon}(\omega)=0.
\]
or, what is the same, for all $\omega\in\Omega_{\theta,2}$
\[
\frac{1}{n}\sum_{k=0}^{n-1} \bigl(
\bigl(\mbox{\upshape Re}\bigl(D_{k}(\theta )e^{ik\theta}\bigr)
\bigr)^{2}I _{[|\mbox{\upshape Re}(D_{k}(\theta)e^{ik\theta})|\geq\epsilon\sqrt
{n}]}+\bigl(\mbox{\upshape Im}\bigl(D
_{k}(\theta)e^{ik\theta}\bigr)\bigr)^{2}I_{[|\mbox{\upshape Im}(D_{k}(\theta
)e^{ik\theta})|
\geq\epsilon\sqrt{n}]}
\bigr)
\]
goes to $0$ in $L^{1}_{\mathbb{P}_{\omega}}$ as $n\to\infty$.

Thus, if $\Omega_{\theta,3}$ is a set of probability one such that
$(\mbox{\upshape Re}(M_{k}(\theta)))_{k\in\mathbb{N}^{*}}$ and
$(\mbox{\upshape Im}(M_{k}(\theta)))_{k
\in\mathbb{N}^{*}}$ is a $(\mathcal{F}_{k-1})_{k\in\mathbb
{N}^{*}}$-adapted martingale in $L^{2}_{\mathbb{P}_{\omega}}$ for all
$\omega\in\Omega_{\theta,3}$ (Lemma~\ref{marher}), the hypotheses 1.
and 2. in Theorem~\ref{bilinvpri} are verified for all $\omega$ in the
set $\Omega_{\theta}$ defined by
%
\begin{equation}
\label{defometheequ} \Omega_{\theta}:=\bigcap_{k=1}^{3}
\Omega_{\theta,k}.
\end{equation}
Since $\mathbb{P}\Omega_{\theta}=1$ this finishes the proof of (a).

To prove (b) we will show that for any given $n\in\mathbb{N}$,
any $\omega\in\Omega_{\theta}$, and any $0\leq t_{1}\leq\cdots
\leq t_{n}$, the $\mathbb{C}^{n}=\mathbb{R}^{2n}$-valued process
\[
\bigl(V_{n}(\theta) (t_{1}), V_{n}(\theta)
(t_{2})-V_{n}(\theta) (t_{1}), \ldots,
V_{n}(\theta) (t_{n})-V_{n}(\theta)
(t_{n-1})\bigr)
\]
has the same asymptotic distribution as
$$\mathbf{B}^{\theta}(t_{1},\ldots, t_{n}):=$$
$$ \bigl[E \bigl\llvert D_{0}(\theta) \bigr\rrvert
^{2}/2 \bigr]^{1/2} \bigl(\Delta_{0}
\mbox{\upshape Re}(B),\Delta_{0}\mbox{\upshape Im}(B),
\Delta_{1}\mbox{\upshape Re}(B), \Delta_{1}
\mbox{\upshape Im}(B), \ldots, \Delta_{n-1}\mbox{\upshape Im}(B) \bigr)$$
under $\mathbb{P}_{\omega}$ where $\Delta_{k}F:=F(t_{k+1})-F(t_{k})$
and $t_{0}=0$ and therefore, by the mapping theorem
(\cite{Bilconpromea}, Theorem~2.7), the finite dimensional asymptotic
distributions of $V_{n}(\theta)$ under $\mathbb{P}_{\omega}$ and those
of (\ref{asyvnthemarcas}) under $\mathbb{P}'$ are the same.

For simplicity, we will assume $n=2$. The argument generalizes easily
to an arbitrary $n\in\mathbb{N}$.

Our goal is thus to prove that for all $\omega\in\Omega_{\theta}$ and
all $0\leq s\leq t$ the asymptotic distribution of
%
\begin{equation}
\label{tesvec} \mathbf{V}_{n}^{\theta}(s,t):=
\bigl(V_{n}(\theta) (s),V_{n}(\theta) (t)-V _{n}(
\theta) (s) \bigr)
\end{equation}
(a $\mathbb{C}^{2}=\mathbb{R}^{4}$-valued process) is the same under
$\mathbb{P}_{\omega}$ as that of
%
$$\mathbf{B}^{\theta}(s,t):=$$
\begin{equation}
\label{asydistesvec}
\bigl[E \bigl\llvert
D_{0}(\theta) \bigr\rrvert^{2}/2 \bigr]^{1/2}\times\bigl( \mbox{\upshape Re}(B) (s),\mbox{\upshape Im}(B) (s), \mbox{\upshape Re}(B)
(t)- \mbox{\upshape Re}(B) (s), \mbox{\upshape Im}(B) (t)-\mbox{\upshape Im}(B) (s)
\bigr)
\end{equation}
under $\mathbb{P}'$.
To prove the convergence in distribution of (\ref{tesvec}) to
(\ref{asydistesvec}), we will use the Cram\'{e}r--Wold theorem. This is,
we will see that for any $\omega\in\Omega_{\theta}$, any
$0\leq s\leq t$, and any
%
\begin{equation}
\label{defu} \mathbf{u}=(a_{1},a_{2},b_{1},b_{2})
\in\mathbb{R}^{4}
\end{equation}
the asymptotic distribution under $\mathbb{P}_{\omega}$ of the
stochastic process $(U_{n})_{n\in\mathbb{N}^{*}}$ defined by
%
\begin{equation}
\label{varcrawoldev} U_{n}:=\mathbf{u}\cdot\mathbf{V}^{\theta}_{n}(s,t)
\end{equation}
is that of a centered normal random variable with variance
%
\begin{equation}
\label{asyvarcrawol} \sigma_{\mathbf{u},s,t}^{2}(\theta):= \frac{E[|D_{0}(\theta)|^{2}]}{2}
\bigl( \bigl(a_{1}^{2}+a_{2}^{2} \bigr)s+
\bigl(b_{1}^{2}+b _{2}^{2} \bigr) (t-s)
\bigr).
\end{equation}

To do so, we will verify the hypotheses of Theorem~\ref{bilmarclt}. Fix
$\mathbf{u}$ as above and note that
\[
U_{n}=\sum_{k=0}^{ \lfloor ns \rfloor}
\eta_{nk}(a_{1},a_{2})+ \sum
_{k= \lfloor ns \rfloor+1}^{ \lfloor nt \rfloor} \eta_{nk}(b_{1},b_{2})
\]
where
%
\begin{equation}
\label{defetank} \eta_{nk}(x_{1},x_{2})=
\frac{1}{\sqrt{n}}(x_{1},x_{2})\cdot e^{ik
\theta}T^{k}D_{0}(
\theta).
\end{equation}

By the construction of $\Omega_{\theta}$, for every $0\leq r$, every
$x_{1}$, $x_{2}$ and every $\omega\in\Omega_{\theta}$, $(\eta_{nk}(x
_{1},\allowbreak x_{2}))_{0\leq k\leq \lfloor nr \rfloor}$ is a triangular
array of $(\mathcal{F}_{k})_{k}$-adapted (real-valued) martingale
differences under $\mathbb{P}_{\omega}$, and by Lemma~\ref{rel16cunmerpel} combined with Remark~\ref{remproone} we can assume
that
%
\begin{equation}
\label{firhyptover} \sum_{k\leq ns}E_{k-1} \bigl[
\eta_{nk}^{2}(a_{1},a_{2}) \bigr]+\sum
_{ns< k\leq nt}E _{k-1} \bigl[\eta_{nk}^{2}(b_{1},b_{2})
\bigr]\to_{n} \sigma_{\mathbf{u},s,t}^{2}( \theta)
\end{equation}
$\mathbb{P}_{\omega}$-a.s.\footnote{More precisely: redefine
$\Omega_{\theta}$ above by intersecting it with the set $\Omega_{
\theta}'$ of elements $\omega$ for which the convergence in Lemma~\ref{rel16cunmerpel} happens $\mathbb{P}_{\omega}$-a.s.} This verifies
the first hypothesis in Theorem~\ref{bilmarclt} under $\mathbb{P}_{
\omega}$ for all $\omega\in\Omega_{\theta}$ for the triangular array
defining $U_{n}$.

It remains to prove that if $\omega\in\Omega_{\theta}$ then
%
\begin{equation}
\label{sechyptover} \sum_{k\leq ns}E_{0} \bigl[
\eta_{nk}^{2}(a_{1},a_{2})I_{[|\eta_{nk}(a_{1},a
_{2})|>\epsilon]}
\bigr](\omega)\to0.
\end{equation}

This is, that for all $\omega\in\Omega_{\theta}$
\[
\sum_{k\leq ns}\eta_{nk}^{2}(a_{1},a_{2})I_{[|\eta_{nk}(a_{1},a_{2})|>
\epsilon]}
\to0
\]
in $L^{1}_{\mathbb{P}_{\omega}}$.

To do so we depart from the Cauchy--Schwarz inequality to get that
\[
\eta_{nk}^{2}(x_{1},x_{2})\leq
\frac
{1}{n}\bigl(x_{1}^{2}+x_{2}^{2}
\bigr)T^{k}\bigl\llvert D _{0}(\theta)\bigr
\rrvert^{2},
\]
so that the sum in (\ref{sechyptover}) is bounded by
\[
\frac{1}{n}\sum_{k\leq ns}E_{0}T^{k}
\bigl[\bigl(a_{1}^{2}+a_{2}^{2}\bigr)\bigl
\llvert D_{0}( \theta)\bigr\rrvert^{2}I_{[(a_{1}^{2}+a_{2}^{2})|D_{0}(\theta)|^{2}\geq
\epsilon^{2} n]}
\bigr].
\]

This obviously goes to zero when $a_{1}=a_{2}=0$. Otherwise it is the
same as
\[
\bigl(a_{1}^{2}+a_{2}^{2}\bigr)
\frac{1}{n}\sum_{k\leq ns}E_{0}T^{k}
\bigl[\bigl\llvert D_{0}(\theta)\bigr\rrvert^{2}I_{[|D_{0}(\theta)|^{2}\geq{\epsilon^{2} n}/{(a_{1}
^{2}+a_{2}^{2})}]}
\bigr],
\]
which, again, goes to zero as $n\to\infty$ for every $\omega\in
\Omega_{\theta}$.
 \qed


\medskip

\begin{cor}
\label{marcasdou}
With the notation introduced in Theorem~\ref{marcas} (see also Remark~\ref{remtwopar}), assume that the function
\[
(\theta,\omega)\mapsto D_{0}(\theta,\omega)
\]
is $\mathcal{B}\otimes\mathcal{F}$-measurable. Then there exists
$\Omega_{0}\subset\Omega$ with $\mathbb{P}\Omega_{0}=1$ such that for
all $\omega_{0}\in\Omega_{0}$ the asymptotic distribution of the
$D[[0,\infty),\mathbb{C}]$-valued function $(\theta,\omega)\mapsto
V_{n}(\theta,\omega)$ corresponds, under $\lambda\times\mathbb{P}
_{\omega_{0}}$, to that of $(\theta,\omega')\mapsto B(\theta,
\omega')$ under $\lambda\times\mathbb{P}'$.
\end{cor}

{\bf Proof:}  The function $(\theta,\omega)\mapsto V_{n}(\theta,
\omega)$ is measurable with respect to $\mathcal{B}\otimes
\mathcal{F}$ (see Section~\ref{thespaD}), and therefore by an
application of Lemma~\ref{mealem} in the appendix and Proposition~\ref{unicon} above, the statement in Theorem~\ref{marcas} can be read
in the following way: for any continuous and bounded function
$f:D[[0,\infty),\mathbb{C}]\to\mathbb{R}$
\[
\lim_{n}E[f\circ V_{n}|\mathcal{B}\otimes
\mathcal{F}_{0}](\theta ,\omega)=E\bigl[f\bigl(B(\theta)\bigr)\bigr]
\]
$\lambda\times\mathbb{P}$-a.s., where the expectation at the right
denotes integration with respect to $\mathbb{P}'$. This is an equality
of $\mathcal{B}\times\mathcal{F}_{0}$ measurable functions, the
$\mathcal{B}$-measurable function at the right being considered as
constant in $\Omega$ for fixed $\theta$. An application of Theorem~34.2(v) in \cite{BilProMea} gives that, for any given $\mathcal{B}
_{0}\subset\mathcal{B}$
%
\begin{equation}
\label{semqueconequ} \lim_{n}E[f\circ V_{n}|
\mathcal{B}_{0}\otimes\mathcal{F}_{0}]=E \bigl[E \bigl[f
\bigl(B( \theta) \bigr) \bigr]|\mathcal{B}_{0}\otimes
\mathcal{F}_{0} \bigr]
\end{equation}
$\lambda\times\mathbb{P}$-a.s. 

If $\mathcal{B}_{0}=\{\emptyset,[0,2
\pi)\}$ is the trivial sigma algebra then, as explained before,
$\lambda_{\theta}:=\lambda$ for all $\theta\in[0,2\pi)$
defines the regular measures corresponding to $E[\,\cdot\,|
\mathcal{B}_{0}]$. It follows, under the light of Corollary~\ref{cormealem}  and Proposition~\ref{unicon}, that
(\ref{semqueconequ}) is nothing but the statement of convergence
$V_{n}\Rightarrow B$ under $\lambda\times\mathbb{P}_{\omega}$ for
$\mathbb{P}$-a.e. $\omega$. This is the desired conclusion.
\qed

\section{Proof of Theorem \protect\ref{Barinvpri}}\label{prothe}
We will present in this and the next section the proofs still pending
from the results announced in Section~\ref{resandcom}. We continue
working under the setting introduced in Section~\ref{genset}.

To begin with, let us prove the following decomposition lemma:

\medskip

\begin{lemma}
\label{firapp}
For all $(n,r,\theta)\in \mathbb{N}\times\mathbb{N}^{*}\times [0,2\pi)$ 
the following equality holds :
\begin{equation}
\label{firappequ}
\begin{array}{lcl}
S_{n}(\theta)-E_{0}S_{n}(\theta)-M_{r,n}(\theta)&=  &-e^{i(n-1)\theta}\left(\sum_{k=1}^{r}  (T^{n-1}E_{0}X_{k}-E_{0}T^{n-1}E_{0}X_{k}) e^{ik\theta}\right)\\
 & & \\
 & +&e^{ir\theta}\sum_{k=2}^{n-1}(T^{k}E_{-1}X_{r}-E_{0}T^{k}E_{-1}X_{r})e^{ik\theta}\\
 & & \\
  & -&D_{r,0}(\theta).
\end{array}
\end{equation}
\end{lemma}

{\bf Proof:} Fix $(n,r,\theta)\in \mathbb{N}\times\mathbb{N}^{*}\times [0,2\pi)$. We depart from the following decomposition of $X_{0}$ (the array is intended to make visible the rearrangements):
\begin{equation}
\label{prolem11}
\begin{array}{rll}
X_{0}=E_{0}X_{0} =&(E_{0}-E_{-1})X_{0} &+ \,\,\,\,\,\,\,\,\,\,\,\,\,\,\,\,\,\,E_{-1}X_{0}\\
 +& (E_{0}-E_{-1})X_{1}e^{i\theta} &- (E_{0}-E_{-1})X_{1}e^{i\theta}+ \\
 +& (E_{0}-E_{-1})X_{2}e^{i2\theta} &- (E_{0}-E_{-1})X_{2}e^{i2\theta}+ \\
  & & \vdots \\
  +& (E_{0}-E_{-1})X_{r}e^{ir\theta} &- (E_{0}-E_{-1})X_{r}e^{ir\theta}\\
  & & \\
=&\sum_{k=0}^{r}(\mathcal{P}_{0}X_{k})e^{ik\theta}&-\sum_{k=1}^{r}(E_{0}X_{k}e^{ik\theta}-E_{-1}X_{k-1}e^{i(k-1)\theta})\\
& & \\
 & & +E_{-1}X_{r}e^{ir\theta}.
\end{array}
\end{equation}
  
Now, using the equality 
$$\sum_{j=0}^{n-1}e^{ij\theta}T^{j}\sum_{k=1}^{r}(E_{0}X_{k}e^{ik\theta}-E_{-1}X_{k-1}e^{i(k-1)\theta})=$$
$$e^{i(n-1)\theta}T^{n-1}\sum_{k=1}^{r}E_{0}X_{k}e^{ik\theta}-\sum_{k=0}^{r-1}E_{-1}X_{k}e^{ik\theta}$$
we get, from (\ref{prolem11}), that
\begin{equation}
\label{prolem12}
\begin{array}{rll}
S_{n}(\theta)= & M_{r,n}(\theta)-( e^{i(n-1)\theta}T^{n-1}\sum_{k=1}^{r}E_{0}X_{k}e^{ik\theta}-\sum_{k=0}^{r-1}E_{-1}X_{k}e^{ik\theta})&
\\

 & & \\
+& \sum_{j=0}^{n-1}e^{ij\theta} T^{j}E_{-1}X_{r}e^{ir\theta}& 
\end{array}
\end{equation}
and that 
\begin{equation}
\label{prolem13}
\begin{array}{rll}
E_{0}S_{n}(\theta)= & D_{r,0}(\theta)-( E_{0}e^{i(n-1)\theta}T^{n-1}\sum_{k=1}^{r}E_{0}X_{k}e^{ik\theta}-\sum_{k=0}^{r-1}E_{-1}X_{k}e^{ik\theta})& \\
 & & \\
+& \sum_{j=0}^{n-1}e^{ij\theta}E_{0}T^{j}E_{-1}X_{r}e^{ir\theta}. &  
\end{array}
\end{equation}

Note that (\ref{firappequ}) follows from (\ref{prolem12}) and (\ref{prolem13}).\qed

\medskip

Let us denote by

\begin{equation}
\label{B}
A_{r,n}=A_{r,n}(\theta,\omega):=\sum_{k=1}^{r}
\bigl(T^{n-1}E_{0}X_{k}(\omega)-E_{0}T^{n-1}E_{0}X
_{k}(\omega) \bigr) e^{ik\theta},
\end{equation}

\begin{equation}
\label{C}
B_{r,n}=B_{r,n}(\theta,\omega):=\sum
_{k=0}^{n-1} \bigl(T^{k}E_{-1}X_{r}(
\omega)-E_{0}T^{k}E_{-1}X_{r}(\omega)
\bigr)e^{ik\theta}.
\end{equation}

The following lemma is a key step for the proof of the validity of our
martingale approximation.

\medskip 

\begin{lemma}
\label{secappequ}
In the context of Lemma~\ref{firapp}, and with the notation (\ref{B})
and (\ref{C}), there exists a constant $C>0$ such that, if $E_{0}$ is
given by the regular version $E_{0}X(\omega)=E^{\omega}X$ ($X\in L
^{1}_{\mathbb{P}}$) then
\begin{enumerate}

\item[1.] For all $(n,r,\theta)\in\mathbb{N}\times\mathbb
{N}^{*}\times[0,2
\pi)$, $\alpha\in\mathbb{R}$, and $\omega\in\Omega$
%
\begin{equation}
\label{firappb} E_{0} \Bigl[ \max_{k\leq n} \bigl
\llvert A_{r,k}(\theta,\cdot) \bigr\rrvert^{2} \Bigr] (
\omega)\leq4\alpha^{2}+ 8\sum_{j=0}^{n-1 }
E_{0}T^{j} \bigl\llvert
\bigl(E_{0}S_{r}(\theta) \bigr)I_{[|E_{0}S_{r}(\theta)|>\alpha]} \bigr
\rrvert ^{2}(\omega).
\end{equation}
\item[2.] For all $\omega\in\Omega$
%
\begin{equation}
\label{firappc}
\int_{0}^{2\pi}{E}_{0} \Bigl[ \max
_{k\leq n} \bigl\llvert B_{r,k}(\theta, \cdot) \bigr
\rrvert^{2} \Bigr] (\omega)\,d\lambda(\theta)\leq C\sum
_{j=2} ^{n-1}E_{{0}}\llvert
E_{j-1}X_{j+r}-E_{0}X_{j+r}\rrvert
^{2}(\omega).
\end{equation}
\end{enumerate}
\end{lemma}

\medskip

{\bf Proof:} We will prove (\ref{firappb}) using a truncation argument:
let $U_{\alpha}$ be the (non-linear) operator given by $U_{\alpha}Y:=
Y I_{|Y|\geq\alpha}$, and let us use the regular version of
$E_{0}$, thus $E_{0}X(\omega)=E^{\omega}X$ ($X\in
L^{1}_{\mathbb{P}}$). Then for all $\omega\in\Omega$

$$\max_{k \leq n} \bigl\llvert A_{r,k}(\theta,\cdot)
\bigr\rrvert(\omega)=\max_{ k \leq n} \bigl\llvert
(Id-E_{0}) \bigl(T^{k-1}E_{0}S_{r}(
\theta) \bigr) \bigr\rrvert^{2}(\omega)\leq$$
$$4\alpha^{2}+ 2\max_{k \leq n} \bigl\llvert
(Id-E_{0})T^{k-1}U_{\alpha} \bigl(E_{0}S
_{r}(\theta) \bigr) \bigr\rrvert^{2}(\omega)\leq$$ $$4 \Biggl(\alpha^{2}+\sum_{j=0}^{n-1}T^{j}
\bigl\llvert U_{\alpha} \bigl(E_{0}S_{r}(\theta)
\bigr) \bigr\rrvert^{2}( \omega)+\sum_{j=0}^{n-1}E_{0}T^{j}
\bigl\llvert U_{\alpha} \bigl(E_{0}S_{r}(\theta)
\bigr) \bigr\rrvert^{2} \Biggr) ( \omega),$$

where we used Jensen's inequality. This clearly implies (\ref{firappb}).

Let us now prove (\ref{firappc}): by Lemma~\ref{hunyou} there exists a
constant $C$ such that
$$\int_{0}^{2\pi}\max_{k\leq n} \bigl
\llvert B_{r,k}(\theta,z) \bigr\rrvert^{2}\,d\lambda( \theta
)\leq C
\int_{0}^{2\pi} \Biggl\llvert\sum
_{j=2}^{n-1} \bigl(T^{j}E_{-1}X_{r}(z)-E
_{0}T^{j}E_{-1}X_{r}(z)
\bigr)e^{ij\theta} \Biggr\rrvert^{2}\,d\lambda(\theta)=$$
$$C\sum_{j=2}^{n-1} \bigl\llvert
E_{j-1}X_{j+r}(z)-E_{0}X_{j+r}(z) \bigr
\rrvert^{2}.$$

The conclusion follows at once by integrating with respect to
$\mathbb{P}_{\omega}$ over these inequalities and using Tonelli's theorem.\qed

\medskip
 
{\bf Proof of Theorem \ref{Barinvpri}:} Consider the notation introduced
in Lemma~\ref{gengenbarinvpridou}. By an application of this lemma and
Corollary~\ref{marcasdou}, it suffices to prove that
%
\begin{equation}
\label{marappp} \tilde{E}_{0} \biggl[ \frac{1}{n}\max
_{1\leq k\leq n} \bigl\llvert S_{k}(\theta, \omega
)-E_{0}S_{k}(\theta,\omega)-M_{r,k}(\theta,
\omega) \bigr\rrvert^{2} \biggr] =o _{r}(1), \qquad\mathbb{P}
\mbox{-a.s.}
\end{equation}

Let us do so: by Lemma~\ref{firapp} and the definition of $\tilde{E}
_{0}$, it is sufficient to prove that there exists $\Omega_{0}$ with
$\mathbb{P}\Omega_{0}=1$ such that if for $(k,r,\theta)\in
\mathbb{N}\times\mathbb{N}^{*}\times[0,2\pi)$ we replace
$Z_{r,k}(\theta,\omega):=A_{r,k}(\theta,\omega)$ or $Z_{r,k}(
\theta,\omega):=B_{r,k}(\theta,\omega)$, then
%
\begin{equation}
\label{apphan} \lim_{r} \limsup_{n}
\int_{0}^{2\pi}{E}_{0} \biggl[
\frac{1}{n} \max_{1\leq k \leq n} \bigl\llvert Z_{r,k}(
\theta,\cdot) \bigr\rrvert^{2} \biggr] (\omega)\,d \lambda(\theta )=0
\end{equation}
for all $\omega\in\Omega_{0}$.

\textit{Proof of (\ref{apphan}) with $Z_{r,k}(\theta,\omega
):=A_{r,k}(\theta,\omega)$}:
if we fix the version of $E_{0}$ given by
$E_{0}X(\omega)=E^{\omega}X$ ($X\in L^{1}_{\mathbb{P}}$) then it is
clear that for any $\omega\in\Omega$
%
\begin{equation}
\label{majrhsb} \bigl\llvert E_{0}S_{r}(\theta
)I_{[|E_{0}S_{r}(\theta)|>\alpha]} \bigr\rrvert(\omega) \leq \Biggl\llvert \Biggl( \sum
_{j=0}^{r-1}E_{{0}}\llvert X_{j}
\rrvert \Biggr) I_{[\sum_{j=0}
^{r-1}E_{0}|X_{j}|>\alpha]} \Biggr\rrvert(\omega),
\end{equation}

and it follows by an application of (\ref{firappb}), Birkhoff's Ergodic
theorem, and Lemma~\ref{dedmerpel} in the {appendix} (fixing first
$\alpha>0$ so that the expectation of the random variable at the right
hand side in (\ref{majrhsb}) is less than any given $\eta>0$) that\vspace*{-2pt}
%
\begin{equation}
\label{btozer} \lim_{n }{E}_{0} \biggl[
\frac{1}{n}\max_{1\leq k \leq n}{|A_{r,k}(\theta,\cdot
)|^{2}} \biggr] =0 \qquad\mathbb{P}\mbox{-a.s.}
\end{equation}
where the (probability one) set $\Omega_{0,1}$ of convergence does not
depend on $\theta$ and, even more, the convergence is uniform in
$\theta$ for any fixed $\omega\in\Omega_{0,1}$. It follows that for
every $\omega\in\Omega_{0,1}$\vspace*{-2pt}

$$\limsup_{n}\int_{0}^{2\pi}{E}_{0}\left[\frac{1}{n}\max_{1\leq k \leq n}{|A_{r,k}(\theta,\cdot)|^{2}}\right](\omega)\,d\lambda(\theta)\leq$$
$$ \int_{0}^{2\pi}\limsup_{n}{E}_{0}\left[\frac{1}{n}\max_{1\leq k \leq n}{|A_{r,k}(\theta,\cdot)|^{2}}\right](\omega)d\lambda(\theta)=0$$
as desired.

\medskip

\textit{Proof of \textup{(\ref{apphan})} with $Z_{r,n}(\theta,\cdot):=B_{r,n}(
\theta,\cdot)$}: again, fix the version of $E_{0}$ given by
$E_{0}X(\omega)=E^{\omega}X$. We depart from (\ref{firappc}) and we note
that, if for every $j\in\mathbb{Z}$, $X_{-\infty,j}:=X_{j}-E_{-
\infty}X_{j}$ then\vspace*{-2pt}
$$\sum_{k=2}^{n-1}E_{0}|(E_{k-1}-E_{0})X_{k+r}|^{2}=\sum_{k=2}^{n-1}E_{0}|(E_{k-1}-E_{0})X_{-\infty,k+r}|^{2}=$$
$$\sum_{k=2}^{n-1}E_{0}T^{k-1}|(E_{0}-E_{-k+1})X_{-\infty,r+1}|^2=\sum_{k=1}^{n-2}(E_{0}T^{k}|E_{0}X_{-\infty,r+1}|^2-|E_{0}X_{-\infty,k+r+1}|^2)\leq$$
$$\sum_{k=1}^{n-2}E_{0}T^{k}|E_{0}X_{-\infty,r+1}|^2 $$
$\mathbb{P}-$a.s. It follows from (\ref{firappc}) and 
Lemma \ref{dedmerpel} that
  \begin{equation}
  \label{ctozer}
  \limsup_{n\to \infty }\int_{0}^{2\pi}{E}_{0}[\frac{1}{n}\max_{1\leq k \leq n}|B_{r,k}(\theta,\cdot)|^2](\omega)\,d\lambda(\theta)\leq C||E_{0}X_{-\infty,r+1}||_{_{\mathbb{P},2}}^2 =C||E_{-(r+1)}X_{-\infty,0}||_{_{\mathbb{P},2}}^2
  \end{equation}
  
$\mathbb{P}-$a.s. over a set $\Omega_{0,2,r}$ independent of $\theta$ and therefore, by the regularity condition (\ref{reg}) (see also (\ref{regconequ}))
  $$\lim_{r}\limsup_{n}\frac{1}{n}\int_{0}^{2\pi}{E}_{0}[\max_{1\leq k \leq n}|B_{r,k}(\theta,\cdot)|^2](\omega)\,d\lambda(\theta)=0$$
for all $\omega\in\Omega_{0,2}:=\bigcap_{r\in\mathbb{N}}\Omega_{0,2,r}$. To conclude this part take $\Omega_{0}:=\Omega_{0,1}\cap\Omega_{0,2}$.\vadjust{\goodbreak}

We are left at this point with the statement on the convergence in
distribution of (\ref{funqueinvprinoncen}) as $n\to\infty$ under the
regularity condition (\ref{reg}).

To prove this claim, start by considering the inequalities ($S_{k}=S
_{k}(\theta,\omega)$ here, and similarly for $M_{r,k}$)
%
$$\tilde{E}_{0}\frac{1}{n}\max_{1\leq k\leq n}\llvert S_{k}-E_{0}S_{k}-M_{r,k}
\rrvert^{2}
\leq2\tilde{E}_{0}\frac{1}{n}\max
_{1\leq k\leq n}\llvert S_{k}-M_{r,k}
\rrvert^{2}+ \frac{2}{n}\tilde{E}_{0}\max
_{1\leq k\leq n}\llvert E_{0}S_{k}\rrvert
^{2}
\leq$$ 
\begin{equation}
\label{ineequcon}
4\tilde{E}_{0}\frac{1}{n}\max_{1\leq k\leq n}
\llvert S_{k}-E_{0}S_{k}-M_{r,k}
\rrvert^{2}+ \frac{6}{n}\tilde{E}_{0}\max
_{1\leq k\leq n}\llvert E_{0}S_{k}\rrvert
^{2}, 
\end{equation}
and notice that, by Lemma~\ref{hunyou} in the {appendix} and $
\mathcal{F}_{0}$-measurability, there exists a constant $C$ such that
\[
\frac{1}{n}\tilde{E}_{0}\max_{1\leq k\leq n}\llvert
E_{0}S_{k}\rrvert ^{2}\leq C \frac{1}{n}\sum
_{k=0}^{n-1}\llvert E_{0}X_{k}
\rrvert^{2},
\]
$\mathbb{P}$-a.s. It follows from these observations and Corollary~\ref{corcarreg} that, under regularity, (\ref{marappp}) is equivalent
to
%
\begin{equation}
\tilde{E}_{0} \biggl[ \frac{1}{n}\max
_{1\leq k\leq n} \bigl\llvert S_{k}(\theta, \omega
)-M_{r,k}(\theta,\omega) \bigr\rrvert^{2} \biggr]
=o_{r}(1), \qquad\mathbb{P}\mbox{-a.s.}
\end{equation}
which gives the convergence of (\ref{funqueinvprinoncen}) under
$\tilde{E}_{0}$ by the argument already given and an easy adaptation of
Lemma~\ref{gengenbarinvpridou}.
\qed

\medskip
\section{Proof of Theorems \protect\ref{barinvprihancon} and \protect\ref{promaxwooae}}\label{prothehancon}
The first step towards the proof of Theorem~\ref{barinvprihancon} is the
following martingale approximation lemma (here $S_{k}$ denotes \textit
{non-rotated} partial sums).

\medskip
\begin{lemma}
\label{marapphanmaxwoo}
In the context of Theorem~\ref{barinvprihancon}
\begin{enumerate}
\item Assuming (\ref{hancon}), there exists $D_{0}\in L^{2}_{\mathbb{P}}(\mathcal{F}_{0})\ominus L^{2}_{\mathbb{P}}(\mathcal{F}_{-1})$ such that, if $M_{n}$ is given by $M_{n}:=\sum_{k=0}^{n-1}T^{k}D_{0}$, then
\begin{equation}
\label{appinehancon}
\lim_{n}E_{0}[\frac{1}{n}\max_{1\leq k\leq n}{|S_{k}-E_{0}S_{k}-M_{k}|^{2}}]=0. \mbox{\,\,\,\,\,\,\,\,\,$\mathbb{P}$\it-a.s.}
\end{equation}
\item Assuming (\ref{maxwoononrot}) there exists $D_{0}\in L^{2}_{\mathbb{P}}(\mathcal{F}_{0})\ominus L^{2}_{\mathbb{P}}(\mathcal{F}_{-1})$ such that, if $M_{n}$ is given by $M_{n}:=\sum_{k=0}^{n-1}T^{k}D_{0}$
\begin{equation}
\label{appinemaxwoo}
\lim_{n}E_{0}[\frac{1}{n}\max_{1\leq k\leq n}{|S_{k}-M_{k}|^{2}}]=0. \mbox{\,\,\,\,\,\,\,\,\,$\mathbb{P}$\it-a.s.}
\end{equation}

\end{enumerate}
\end{lemma}

{\bf Proof:} The first statement is a part of Theorem  2.3 in \cite{CuVo}. The second statement is Theorem 2.7 in \cite{cunmer}.

\medskip

\begin{remark}
\label{remmarappequundmaxwoo}
If $(a_{n})_{n}$, $(b_{n})_{n}$ are sequences of complex  numbers with $|b_{n}|$ increasing to $\infty$  and $a_{n}/b_{n}\to 0$, then 
$$\lim_{n}\max_{1\leq k\leq n}\frac{|a_{k}|}{|b_{n}|}=0$$
because 
$$0\leq\limsup_{n}\max_{1\leq k\leq n}\frac{|a_{k}|}{|b_{n}|}\leq \limsup_{N}\limsup_{n}(\frac{1}{|b_{n}|}\sum_{k=1}^{N-1}{|a_{k}|}+ \sup_{k\geq N}\frac{|a_{k}|}{|b_{k}|})=$$
$$\lim_{N}\sup_{k\geq N}\frac{|a_{k}|}{|b_{k}|}=0.$$
If we apply this observation to the inequalities (\ref{ineequcon}) with
$\tilde{E}_{0}$ replaced by $E_{0}$ we see that, under the hypothesis
%
\begin{equation}
\label{contozerhyp} \lim_{n}\frac{E_{0}S_{n}}{\sqrt{n}}=0 \qquad
\mathbb{P}\mbox{-a.s.},
\end{equation}
(\ref{appinehancon}) and (\ref{appinemaxwoo}) are equivalent. We can
actually verify (\ref{contozerhyp}) under (\ref{maxwoononrot}) (see
\cite{cunmer}, Proposition~4.9) and therefore, \textit{under the
(non-rotated) Maxwell and Woodroofe condition} (\ref{maxwoononrot}),
(\ref{appinehancon}) \textit{and} (\ref{appinemaxwoo}) \textit{imply each other}.
\end{remark}

\medskip

\begin{remark}
\label{remmarappmaxwoohan}
It is convenient, for the sake of clarity, to describe explicitly the martingale differences stated in Lemma \ref{marapphanmaxwoo}: under (\ref{maxwoononrot}), $D_{0}$ can be described as the $\mathbb{P}-$a.s. and in $L^{2}_{\mathbb{P}}-$convergent series
\begin{equation}
\label{mardifmaxwoo}
D_{0}:=\sum_{n\geq 0} \sum_{k\geq n}\frac{\mathcal{P}_{0}X_{k}}{k+1},
\end{equation}
and under (\ref{hancon}) $D_{0}$ is given by the $\mathbb{P}-$a.s. and in $L^{2}_{\mathbb{P}}-$convergent series
\begin{equation}
\label{mardifhancon}
D_{0}:=\sum_{k\geq 0}\mathcal{P}_{0}X_{k}.
\end{equation}
\end{remark} 

\medskip

{\bf Proof of Theorem \ref{barinvprihancon}\,:} Start by considering the following observations: let $(\Omega,\mathcal{F},\mathbb{P})$ be the underlying probability  space  (the domain of $X_{0}$), and for every measurable function  $Y:\Omega\to\mathbb{C}$
let $\tilde{Y}$ be the {\it extension to the product space} $\tilde{Y}:[0,2\pi)\times\Omega\to\mathbb{C}$ defined by 
$$\tilde{Y}(u,\omega):=e^{iu}Y(\omega).$$

For $\theta\in [0,2\pi)$ fixed, define $\tilde{T}_{\theta}:[0,2\pi)\times\Omega\to [0,2\pi)\times\Omega$ by $$\tilde{T}_{\theta}(u,\omega):=((u+\theta)mod(2\pi),T\omega)$$ 
($\tilde{T}_{\theta}$ is the product map between the rotation by an angle of $\theta$ and $T$). It is easy to see that $\tilde{T}_{\theta}$ is invertible, that it preserves the product measure $\lambda\times\mathbb{P}$, and that it is $\mathcal{B}\otimes\mathcal{F}-$bimeasurable. It is also easy to see that for every (measurable) $Y:\Omega\to\mathbb{C}$ and every $k\in\mathbb{Z}$, $\tilde{T}_{\theta}^{k}\tilde{Y}=\widetilde{e^{ik\theta}T^{k}Y}$, and that $||\tilde{Y}||_{_{\lambda\times{\mathbb{P}},2}}=||Y||_{_{\mathbb{P},2}}$ if $Y\in L^{2}_{\mathbb{P}}$.

Denote, for every $k\in\mathbb{Z}$, $\tilde{\mathcal{F}}_{k}:=
\mathcal{B}\otimes\mathcal{F}_{k}$, $\tilde{E}_{k}:={E}[\cdot|
\tilde{\mathcal{F}}_{k}]$ (conditional expectation with respect to
$\lambda\times\mathbb{P}$), and $\tilde{\mathcal{P}}_{k}:=\tilde{E}
_{k}-\tilde{E}_{k-1}$ (avoid confusion with the meaning of this notation
in Lemma~\ref{gengenbarinvpridou}: here $\mathcal{B}_{0}=\mathcal{B}$).
It is easy to see that $(\tilde{\mathcal{F}}_{k})_{k\in\mathbb{Z}}$ is
a $\tilde{T}_{\theta}$-filtration.

Recall also the notation introduced in Definition~\ref{defdisfoutra},
and notice now that, if $\theta\in(0,2\pi)$ and $Y\in L^{1}_{
\mathbb{P}}$ is given then
%
\begin{equation}
\label{froimphantohan} \bigl(1-e^{i\theta} \bigr)S_{n}(Y,T,\theta,\cdot
)=T^{-1}Y+S_{n} \bigl(Y-T^{-1}Y,T, \theta,\cdot
\bigr)-T^{n-1}Y e^{in\theta}.
\end{equation}

Theorem~\ref{barinvprihancon} consists of two statements: the
convergence of (\ref{funqueinvpri}) under (\ref{imhancon}) or
(\ref{maxwoo}) and the convergence of (\ref{funqueinvprinoncen}) under
(\ref{maxwoo}). For reasons inherent to the logic of our arguments, we
will prove first the convergence of (\ref{funqueinvpri}) under
(\ref{imhancon}), and then the convergence of (\ref{funqueinvpri})
\textit{and} (\ref{funqueinvprinoncen}) under (\ref{maxwoo}).

\textit{Convergence of} (\ref{funqueinvpri}) \textit{under}
(\ref{imhancon}). First note that, \textit{without loss of generality, we
can assume that $(X_{k})_{k\in\mathbb{Z}}$ is regular (Definition
\textup{\ref{regcondef}})}. Indeed, with the notation introduced in Definition~\ref{defdisfoutra}:
\begin{equation}
\label{redregcas}
S_{k}(X_{0}, T,\theta,\cdot)-E_{0}S_{k}(X_{0}, T,\theta,\cdot)=S_{k}(X_{0}-E_{-\infty}X_{0}, T,\theta,\cdot)-E_{0}S_{k}(X_{0}-E_{-\infty}X_{0}, T,\theta,\cdot),
\end{equation}
and therefore we can study the desired asymptotics replacing the stationary process $(X_{k})_{k\in\mathbb{Z}}$ by the (stationary and) regular process $(X_{k}-E_{-\infty}X_{k})_{k\in\mathbb{Z}}$. We will therefore assume in this part of the proof that $(X_{k})_{k\in\mathbb{Z}}$ is, indeed, regular.

\textit{Martingale approximations.} For the sake of clarity, let us depart
from the following observation: if we assume that
\[
\Delta_{0}(\theta):=\sum_{k\geq0}
\mathcal{P}_{0}(X_{k}-X_{k-1})e ^{ik\theta}
\]
converges $\mathbb{P}$-a.s. and in $L^{2}_{\mathbb{P}}$ then,
necessarily
%
\begin{equation}
\label{mardifprogenhan} D_{0}(\theta):=\sum_{k\geq0}
\mathcal{P}_{0}X_{k}e^{ik\theta}
\end{equation}
converges in $L^{2}_{\mathbb{P}}$ (under regularity). Even more,
%
\begin{equation}
\label{reldelzerdzer} \bigl(1-e^{i\theta} \bigr)D_{0}(\theta
)=\Delta_{0}(\theta),
\end{equation}
$\mathbb{P}$-a.s. To see this apply (\ref{froimphantohan}) and the definition of $\mathcal{P}_{0}$ to obtain,
for every $r\in\mathbb{N}$, that
\[
\bigl(1-e ^{i\theta}\bigr)\sum_{k=0}^{r}\mathcal{P}_{0}X_{k}e^{ik\theta}=\sum_{k=0}^{r}
\mathcal{P}_{0}(X_{k}-X_{k-1})e^{ik\theta}-
\mathcal{P}_{0}X_{r}e^{i(r+1)\theta},
\]
and let $r\to\infty$ taking into account the regularity of
$(X_{k})_{k\in\mathbb{Z}}$.

Now note that (\ref{imhancon}) implies that 
\begin{equation}
\label{imhanconlif}
\sum_{k\geq 0}||{\mathcal{P}}_{0}({X}_{k}-{X}_{k-1})||_{_{\mathbb{P},2}}=\sum_{k\geq 0}||\tilde{\mathcal{P}}_{0}(\tilde{X}_{k}-\tilde{X}_{k-1})||_{_{\lambda\times\mathbb{P},2}}<\infty.
\end{equation}

Taking $Y_{0}=X_{0}-X_{-1}$ and noticing that, 
$$||\tilde{\mathcal{P}}_{0}\tilde{T}^{k}_{\theta}\tilde{Y}_{0}||_{_{\lambda\times\mathbb{P},2}}=||\mathcal{P}_{0}T^{k}Y_{0}||_{_{\mathbb{P},2}}= ||{\mathcal{P}}_{0}({X}_{k}-{X}_{k-1})||_{_{\mathbb{P},2}}$$ 
we see by Lemma \ref{marapphanmaxwoo} that we can find $\tilde{\Delta}_{0}\in L^{2}_{\lambda\times\mathbb{P}}(\tilde{\mathcal{F}}_{0})\ominus L^{2}_{\lambda\times\mathbb{P}}(\tilde{\mathcal{F}}_{-1})$ such that
\begin{equation}
\label{appimhanconlifequ}
\lim_{n}\frac{1}{n}\tilde{E}_{0}\max_{1\leq k\leq n}|S_{k}(\tilde{X}_{0}-\tilde{X}_{-1},\tilde{T}_{\theta},0,\cdot)-E_{0}S_{k}(\tilde{X}_{0}-\tilde{X}_{-1},\tilde{T}_{\theta},0,\cdot)-S_{k}(\tilde{\Delta}_{0},\tilde{T}_{\theta},0,\cdot)|^2=0,
\end{equation}
$\lambda\times\mathbb{P}-$a.s. Even more, using Remark \ref{remmarappmaxwoohan} we see that 
\begin{equation}
\label{mardifimhanconlifequ}
\tilde{\Delta}_{0}(u,\omega)=e^{iu}\sum_{k\geq 0}{\mathcal{P}_{0}(X_{k}-X_{k-1})(\omega)e^{ik\theta}}.
\end{equation} 
This gives that if $\Delta_{0}(\theta,\omega):=e^{-iu}\tilde{\Delta}_{0}(u,\omega)$ then
\begin{equation}
\label{appimhanconequ}
\lim_{n}\frac{1}{n}{E}_{0}\max_{1\leq k\leq n}|S_{k}({X}_{0}-{X}_{-1},{T},\theta,\cdot)-E_{0}S_{k}({X}_{0}-{X}_{-1},{T},\theta,\cdot)-S_{k}({\Delta}_{0}(\theta),{T},\theta,\cdot)|^{2}=0,
\end{equation}
$\mathbb{P}$-a.s., which is the same, by (\ref{froimphantohan}), (\ref{reldelzerdzer}) and the convergence $E_{0}|(Id-E_{0})X_{n}|^2=o(n)$, $\mathbb{P}$-a.s. (use for instance the pointwise ergodic theorem for $E_{0}T$, see also the first part of  Remark~\ref{remmarappequundmaxwoo}),  as
%
\begin{eqnarray}
\label{appimhanconequalmost}
&&\bigl|1-e^{i\theta} \bigr|^{2}\lim_{n}
\frac{1}{n}{E}_{0}\max_{1\leq k\leq n} \bigl\llvert
\bigl(S_{k}( {X}_{0},{T},\theta,\cdot)-E_{0}S_{k}({X}_{0},{T},
\theta,\cdot)-S _{k} \bigl({D}_{0}(\theta),{T},\theta, \cdot
\bigr) \bigr) \bigr\rrvert^{2}\nonumber\\[-8pt]\\[-8pt]\nonumber
&&\quad =0,
\end{eqnarray}
$\mathbb{P}$-a.s., where $D_{0}(\theta)$ is given by
(\ref{mardifprogenhan}). The result follows at once from Lemma~\ref{gengenbarinvpri} assuming that $e^{2i\theta}\notin \mathrm{Spec}_{p}(T)$
(in particular $e^{i\theta}\neq1$).

\textit{Convergence of \textup{(\ref{funqueinvpri})} and
\textup{(\ref{funqueinvprinoncen})} under \textup{(\ref{maxwoo})}}. With the notation
already introduced, note that
%
\begin{equation}
\label{maxwoolif} \sum_{k\geq1}\frac{\|
{E}_{0}{S}_{k}({X}_{0},{T},\theta,\cdot)\|_{
\mathbb{P},2}}{k^{3/2}}=\sum
_{k\geq1}\frac{\|\tilde{E}_{0}{S}_{k}(
\tilde{X}_{0},\tilde{T}_{\theta}, 0,\cdot)\|_{\lambda\times
\mathbb{P},2}}{k^{3/2}}<\infty
\end{equation}
and it follows, by applying Lemma~\ref{marapphanmaxwoo} again, that if
$\tilde{D}_{0}$ is defined by
%
\begin{equation}
\label{mardifmaxwoolifequ} \tilde{D}_{0}:=\sum_{n\geq0}
\sum_{k\geq n}\frac{\tilde{\mathcal{P}}
_{0}\tilde{T}^{k}_{\theta}\tilde{X}_{0}}{k+1},
\end{equation}
then
%
\begin{equation}
\label{appmaxwoolifequ} \lim_{n}\frac{1}{n}
\tilde{E}_{0}\max_{1\leq k\leq n} \bigl\llvert
S_{k}( \tilde{X} _{0},\tilde{T}_{\theta},0,
\cdot)-S_{k}( \tilde{D}_{0},\tilde{T}_{
\theta},0,
\cdot) \bigr\rrvert ^{2}=0
\end{equation}
$\lambda\times\mathbb{P}$-a.s. By arguments similar to those preceding
(\ref{appimhanconequalmost}), this is the same as saying that if
%
\begin{equation}
\label{mardifmaxwoofixfre} {D}_{0}(\theta)=\sum_{n\geq0}
\sum_{k\geq n}\frac{\mathcal{P}_{0}X
_{k}e^{ik\theta}}{k+1}
\end{equation}
then
%
\begin{equation}
\label{appmaxwooequ} \lim_{n}\frac{1}{n}{E}_{0}
\max_{1\leq k\leq n} \bigl\llvert S_{k}({X}_{0},{T},
{\theta},\cdot)-S_{k} \bigl({D}_{0}(\theta),{T},{\theta},
\cdot \bigr) \bigr\rrvert^{2}=0,
\end{equation}
$\mathbb{P}$-a.s. This implies the desired conclusion by (an easy
adaptation of) Lemma~\ref{gengenbarinvpri}. The convergence of
(\ref{funqueinvpri}) follows by a similar argument after applying Remark~\ref{remmarappequundmaxwoo} to (\ref{appmaxwoolifequ}).
\qed

{\bf Proof of Theorem~\ref{promaxwooae}:} We start from the
following observation: if $(a_{k})_{k\in\mathbb{N}}$ is a sequence of
nonnegative numbers, $(b_{k})_{k\in\mathbb{N}}$ is a sequence of
(strictly) positive numbers, and $n\in\mathbb{N}^{*}$ is given, an
application of H\"{o}lder's inequality gives that
\[
\sum_{k=1}^{n}\frac{a_{k}}{k^{3/2}}\leq
\Biggl( \sum_{k=1}^{n}\frac{1}{kb
_{k}}
\Biggr) ^{1/2} \Biggl( \sum_{k=1}^{n}
\frac{a_{k}^{2}b_{k}}{k^{2}} \Biggr) ^{1/2},
\]
and therefore the convergence of $\sum_{k}a_{k}/k^{3/2}$ ($k\geq1$) is
equivalent to the existence of a positive sequence $(b_{k})_{k\in
\mathbb{N}}$ such that both $\sum_{k}1/(kb_{k})$ and $\sum_{k}a_{k}
^{2}b_{k}/k^{2}$ are convergent (for the necessity consider
$b_{k}:=k^{1/2}/a_{k}$ if $a_{k}>0$ and $b_{k}=k$ otherwise).

This observation applied to $a_{k}:=\|E_{0}S_{k}(\theta)\|_{
\mathbb{P},2}$ (for $\theta\in[0,2\pi)$ fixed) gives that a
sufficient condition for the fulfillment of (\ref{maxwoo}) is the
existence of $\beta>1$ such that
%
\begin{equation}
\label{appmaxwooint} \sum_{k\in\mathbb{N}^{*}}(\log k)^{\beta}
\frac{E|E_{0}S_{k}(\theta
)|^{2}}{k^{2}}<\infty.
\end{equation}
In order to provide a condition giving rise to the fulfillment of
(\ref{appmaxwooint}) for $\lambda$-a.e. $\theta$, we fix $n\in
\mathbb{N}^{*}$ and start by noticing that, by orthogonality and
Fubini's theorem
$$
\int\sum_{k=1}^{n}(\log
k)^{\beta}\frac{\|E_{0}S_{k}(\theta)\|_{
\mathbb{P},2}^{2}}{k^{2}}\,d\lambda(\theta)=\sum
_{k=1}^{n}\sum_{j=0}^{k-1}
\frac{(\log k)^{\beta}}{k^{2}}\|E_{0}X_{j}\|_{
\mathbb{P},2}^{2}$$
$$= \sum_{j=0}^{n-1}||E_{0}X_{j}||_{\mathbb{P},2}^{2}
\sum_{k=j+1}^{n}\frac{(
\log k)^{\beta}}{k^{2}},
$$
and that one can show (using for instance the integral test and
integration by parts) that for every $\beta>1$, there exists
$C(\beta)>0$ such that
\[
\sum_{k=j+1}^{n}\frac{(\log k)^{\beta}}{k^{2}}\leq C(
\beta)\frac{(
\log(j))^{\beta}}{j},
\]
for every $1\leq j\leq n$ given. This, together with Tonelli's theorem,
shows that
\[
\int\sum_{k=1}^{\infty}(\log
k)^{\beta}\frac{\|E_{0}S_{k}(\theta)\|_{
\mathbb{P},2}^{2}}{k^{2}}\,d\lambda(\theta)\leq C(\beta)\sum
_{k=1}^{\infty}(\log k)^{\beta}\frac{\|E_{0}X_{k}\|_{\mathbb{P},2}
^{2}}{k},
\]
$\mathbb{P}$-a.s. With this Theorem~\ref{promaxwooae} follows from
(\ref{equpromaxwooae}) and Theorem~\ref{barinvprihancon}.\qed

\renewcommand\thesection{A}
\section*{Appendix}
\label{appp}
\renewcommand{\theequation}{A.\arabic{equation}}
\setcounter{equation}{0}
In this section, we provide some results used along the proofs of the
statements previously given. Some of these results belong to the
existing literature and are included here for the sake of clarity, the
rest of them are either not very visible in the mainstream literature or new,
and we include them in this section due to their general scope.

We start by giving a further equivalence to the Portmanteau theorem,
valid in the case of separable metric spaces, and whose relevance for
our arguments lies in the fact that it reduces the ``integral testing'' for convergence in distribution
to a countable set of functions.

To introduce this result, first remember the notion of a {Urysohn
function}: given two closed, disjoint sets $F_{0}$, $F_{1}$ in a perfectly
normal topological space (for instance, any metric space) $
\mathcal{T}$,
\[
U=U(F_{0},F_{1}):\mathcal{T}\to[0,1]
\]
is called a \textit{Urysohn function} if it is continuous,
$U^{-1}\{0\}=F_{0}$ and $U^{-1}\{1\}=F_{1}$.

Let us call a collection $\{F_{j}\}_{j\in J}$ of closed sets in
$\mathcal{T}$ a \textit{co-base} if $\{\mathcal{T}\setminus
F_{j}\}_{j\in J}$ is a base of $\mathcal{T}$. We will also use the following notation: if $S$ is a metric space with
distance function $d$, then for any given $x\in S$ and $A\subset S$ (not
necessarily in the topology of $S$) we define the \textit{distance from $x$ to $A$} by
\[
d(x,A):=\inf_{a\in A}\,d(x,a),
\]
and we define the \textit{$\epsilon$-neighborhood of
$A$}, $A^{\epsilon}$, as the (open) set
\[
A^{\epsilon}:=\bigl\{x\in S:d(x,A)<\epsilon\bigr\}.
\]

\begin{lemmaapp}
\label{imppro}
Let $S$ be a separable metric space. Denote by $\mathbf{C}^{b}(S)$ the
space of bounded, continuous real-valued functions on $S$. Let
$\{F_{n}\}_{n\in\mathbb{N}}$ be a co-base of $S$ which is also a
$\pi$-system, and let $X_{n}$, $X$ ($n\in\mathbb{N}$) be random
elements of $S$ (not necessarily defined on the same probability space).
Then the following two statements are equivalent
\begin{enumerate}
\item For every $f\in\mathbf{C}^{b}(S)$,
\[
\lim_{n} Ef(X_{n})=Ef(X).
\]
\item For every $k\in\mathbb{N}$, every rational $\epsilon>0$,
and some
Urysohn function $U_{k,\epsilon}=U(S\setminus F_{k}^{\epsilon},F
_{k})$
\[
\lim_{n}EU_{k,\epsilon}(X_{n})=EU_{k,\epsilon}(X).
\]
\end{enumerate}
\end{lemmaapp}

{\bf Proof:} Denote by $P_{n}$ the law of $X_{n}$ and by $P$ the law of $X$. Since 1. clearly implies 2. it suffices to see, by the Portmanteau Theorem (\cite{Bilconpromea}, Theorem 2.1), that if 2. is true then for any given closed set $F$
$$\limsup_{n}{P}_{n}F\leq PF.$$
If for some $k$, $F=F_{k}$, this is a consequence of the inequalities
$$I_{F}\leq U_{k,\epsilon}\leq I_{F^{\epsilon}},$$
the hypothesis in 2. and the continuity from above of finite measures. 

If $F$ is an arbitrary closed set, say $F=\cap_{j\in J}F_{j}$ for some $J\subset \mathbb{N}$, and if we define for all $k\in \mathbb{N}$, $J_{k}:=J\cap [0,k]$ and $A_{k}:=\cap_{j\in J_{k}} F_{j}$ then, since $A_{k}\in\{F_{n}\}_{n}$,
$$\limsup_{n} P_{n}F\leq \limsup_{n} P_{n} A_{k} \leq PA_{k} $$ 
for all $k$. By letting $k\to \infty$ we get the desired conclusion.\qed

We remark that the Portmanteau theorem can be extended to the context of abstract perfectly normal spaces if one inteprets convergence in distribution as the fulfillment of the hypothesis 1. of Lemma \ref{imppro}.  This can be seen by following the arguments in \cite{Bilconpromea} and using the fact that every closed set is a $G_{\delta}$ set. In this context Lemma \ref{imppro} corresponds to the second-countable case.

\medskip

Our next result, Theorem \ref{limlimlem}, is an improvement due to Dehling, Durieu and Volny, of Theorem 3.1 in \cite{Bilconpromea} for the case in which the target (state) space is a complete and separable metric space. As in the previous pages, ``$\Rightarrow$'' denotes convergence in distribution here.

\medskip

\begin{thmapp}
\label{limlimlem}
{\it Let $(S,{d})$ be a complete and separable metric space. Assume that for all $r,n\geq 0$, $X_{(r,n)}$ and $X_{n}$ are random elements of $S$ defined on the same probability space $(\Omega,\mathcal{F},\mu)$, and that $X_{(r,n)}\Rightarrow_{n} Z_{r}$. Then the hypothesis
\begin{equation}
\label{liminf} \lim_{r}\limsup_{n}\mu
\bigl[d(X_{(r,n)},X_{n})\geq\epsilon \bigr]=0\quad\quad
\textit{for all }\epsilon>0
\end{equation}
implies the existence of a random element $X$ of $S$ such that $Z_{r}\Rightarrow_{r} X$  and  $X_{n}\Rightarrow_{n}X$.}
\end{thmapp}

{\bf Proof:} This is Theorem 2 in \cite{deduvo}.\qed
\medskip

\begin{corapp}
\label{corlimlimlem}
In the context of Theorem \ref{limlimlem} denote, for any given $q>0$, $$||Z||_{\mu,q}:=\left(\int_{\Omega}|Z|^{q}d\mu(\omega)\right)^{1/q}.$$ 
If for some $q> 0$
$$\lim_{r}\limsup_{n}||d(X_{(r,n)},X_{n})||_{\mu,q}=0$$
and $X_{(r,n)}\Rightarrow_{n} Z_{r}$, then there exists a random element $X$ such that $X_{n}\Rightarrow_{n} X$.
\end{corapp}

{\bf Proof:} Apply Markov's inequality to verify the hypothesis of Theorem \ref{limlimlem}. \qed

\medskip

The following lemma is used without a proof along the references consulted by the author, thus a demonstration is given. This result allows us to pass from the study of stationary martingales to martingales under the conditional regular measures (see Section \ref{genset}).

\medskip

\begin{lemmaapp}
\label{marher}
With the notation and definitions given in Section~\ref{genset}, and
denoting further by $E^{\omega}_{k}$ the conditional expectation with
respect to $\mathcal{F}_{k}$ and $\mathbb{P}_{\omega}$, the following
property holds: for every $k\in\mathbb{N}$,  every
$\mathbb{P}$-integrable $Y$, and every fixed version of $E_{k}Y$ (also denoted by $E_{k}Y$):
%
\begin{equation}
\label{herconexp} E^{\omega}_{k}Y=E_{k}Y
\end{equation}
$\mathbb{P}_{\omega}$-a.s. for $\mathbb{P}$-a.e. $\omega$. In
particular, if $(M_{n})_{n}$ is an $(\mathcal{F}_{n})_{n}$-adapted
martingale in $L^{p}_{\mathbb{P}}$ ($p\geq1$), then $(M_{n})_{n
\geq0}$ is an $(\mathcal{F}_{n})_{n}$-adapted martingale in
$L^{p}_{\mathbb{P}_{\omega}}$ for $\mathbb{P}$-a.e.~$\omega$.
\end{lemmaapp}

\medskip

{\bf Proof:} Fix a version of $Y\in L^{1}_{\mathbb{P}}$. We will prove that for any ($\mathcal{F}_{k}-$measurable) version of $E_{k}Y$, there exists a set $\Omega_{Y}\subset \Omega$ with $\mathbb{P}\Omega_{Y}=1$  such that the following holds: for every $\omega \in \Omega_{Y}$ and every $A\in \mathcal{F}_{k}$
\begin{equation}
\label{equconexp}
\int_{A} Y(z)d\mathbb{P}_{\omega}(z)=\int_{A} E_{k}Y(z)d\mathbb{P}_{\omega}(z),
\end{equation}
this clearly implies the first conclusion.

Fix a ($\mathcal{F}_{k}-$measurable) version of $E_{k}Y$ and notice that for $A$ fixed, a set $\Omega_{Y,A}$ of probability one such that (\ref{equconexp}) holds for all $\omega\in\Omega_{Y,A}$ exists by the property defining the family $\{\mathbb{P}_{\omega}\}_{\omega\in \Omega}$ and because
$$E_{0}[YI_{A}]=E_{0}[(E_{k}Y)I_{A}],$$
$\mathbb{P}-$a.s. Without loss of generality $\Omega_{Y,A}\subset \{\omega\in \Omega: |Y|+|E_{k}Y|\in L^{1}_{\mathbb{P}_{\omega}}\}$ (the last set has $\mathbb{P}-$measure one because $E|Z|=EE_{0}|Z|$ for every $Z\in L^{1}_{\mathbb{P}}$).

Now proceed as follows: let $\{A_{n}\}_{n\in\mathbb{N}}\subset \mathcal{F}_{k}$ be a countable family generating $\mathcal{F}_{k}$  which is also a $\pi-$system and includes $\Omega$ (such a family exists because $\mathcal{F}_{0}$ is assumed countably generated), let $\Omega_{Y}:=\cap_{n\geq 1}\Omega_{Y,A_{n}}$, and let $\mathcal{G}_{k}\subset\mathcal{F}_{k}$ be the family of sets $A \in \mathcal{F}_{k}$ such that ($\ref{equconexp}$) holds for all $\omega \in \Omega_{Y}$. It is easy to see that $\mathcal{G}_{k}$ is a $\lambda-$system and therefore, since it includes $\{A_{n}\}_{n\in\mathbb{N}}$, the $\pi-\lambda$ theorem implies that $\mathcal{G}_{k}=\mathcal{F}_{k}$. Note that $\mathbb{P}\Omega_{Y}=1$, and that for all $\omega\in \Omega_{Y}$, (\ref{equconexp}) holds for all $A\in \mathcal{F}_{k}$. 

This gives the proof of the first conclusion. 
The second conclusion (the one about martingales) follows easily from this, together with the fact that  $E|X|^{p}=EE_{0}|X|^{p}$ and therefore $E|X|^{p}<\infty$ if and only if $E^{\omega}|X|^{p}<\infty$ for $\mathbb{P}-$a.e. $\omega$.\qed

\medskip

Recall the following (Doob's) maximal inequality (\cite{RevYor}, p.53): {if $p>1$ is given and $(M_{n})_{n\geq 0}$ is a positive submartingale in $L^{p}_{\mu}$ then}
\begin{equation}
\label{doo} \llVert M_{n}\rrVert_{p,\mu}\leq \Bigl
\llVert\max_{0\leq k\leq n}M_{k} \Bigr\rrVert_{p,\mu}
\leq\frac{p}{p-1}\llVert M_{n}\rrVert_{p,\mu}.
\end{equation}
A combination of Doob's maximal inequality (\ref{doo}) with Lemma \ref{marher} gives the following lemma.

\medskip

\begin{lemmaapp}
\label{doomaxl2que}
With the notation of Section~\ref{genset}, if $(M_{k})_{k}$ is an
$(\mathcal{F}_{k})_{k}$-adapted complex martingale in $L^{2}_{
\mathbb{P}}$ then
%
\begin{equation}
\label{Doo} E_{0} \Bigl[\max_{0\leq k\leq n}\llvert
M_{k}\rrvert \Bigr]^{2} \leq4E_{0}\llvert
M_{n}\rrvert^{2}, \qquad\mathbb{P}\mbox{-a.s.}
\end{equation}
\end{lemmaapp}

We also need in this paper the following ergodic theorem, which was demonstrated at the beginning of the proof of Propostion \ref{proconsqu}.

\medskip
\begin{lemmaapp}
\label{dedmerpel}
With the notation introduced in Section~\ref{genset}, and assuming
$T$ is ergodic, for every $Y\in L^{1}_{\mathbb{P}}$,
%
\begin{equation}
\label{ergtheequ} \lim_{n}\frac{1}{n}\sum
_{k=0}^{n-1}E_{0}T^{k}Y= EY,
\end{equation}
$\mathbb{P}$-a.s. and in $L^{1}_{\mathbb{P}}$.
\end{lemmaapp}

{\bf Proof:} 
See the proof of Proposition~\ref{proconsqu}.\hskip.2pt\footnote{Or see Lemma~7.1 in  \cite{DedMerPel} for a slightly different version of this result (which also inspired it).}\qed

\bigskip

The following lemma is a classical tool in Harmonic Analysis, we give
here a concrete version sufficient for our purposes.

\medskip

\begin{lemmaapp}
\label{hunyou}
There exists a constant $C$ with the following property: for any given
$f\in L^{2}_{\lambda}$ with Fourier expansion
\[
S_{f}(\theta)=\sum_{k\geq0}a_{k}e^{ik\theta},
\]
and denoting by $S_{f,n}(\theta):=\sum_{k=0}^{n-1}a_{k}e^{ik\theta}$
the $n$th Fourier partial Fourier sum of $f$:
\[
\int_{[0,2\pi)}\sup_{n}\bigl\llvert
S_{f,n}(\theta)\bigr\rrvert^{2} \,d\lambda (\theta) \leq C
\int_{[0,2\pi)}\bigl\llvert f(\theta)\bigr\rrvert^{2} \,d
\lambda (\theta).
\]
\end{lemmaapp}

{\bf Proof:} \cite{hunyou}.\qed 

\bigskip

The next lemma is useful to compute the finite-dimensional asymptotic
distributions that identify our asymptotic (quenched) limits.

\medskip
\begin{lemmaapp}
\label{lemcunmerpel}
Let $(\Omega,\mathcal{F},\mathbb{P})$ be a probability space, let
$T:\Omega\to\Omega$ be a measure preserving transformation and let
$\theta\in\mathbb{R}$. If the only integrable (complex-valued)
function $Y$ satisfying $TY=e^{-i\theta}Y$ is $Y=0$ (i.e., if $e^{-i\theta}\notin \mathrm{Spec}_{p}(T)$), then for every
$X\in L^{1}_{\mathbb{P}}(\Omega,\mathbb{C})$
%
\begin{equation}
\label{ergthedisfoutraequ} \lim_{n}\frac{1}{n}\sum
_{k=0}^{n-1}T^{k}Xe^{ik\theta}=0,\qquad
\mathbb{P}\mbox{-a.s. and in }L^{1}_{\mathbb{P}}.
\end{equation}
\end{lemmaapp}

{\bf Proof:} \cite{CuMePe}, p.~20.\footnote{Actually, as proved for
instance in \cite{bardis}, Theorem~3.2, if $X\in L^{p}_{\mathbb{P}}$
for some $p\geq1$ and $\theta$ is arbitrary, the random variables at
the left in (\ref{ergthedisfoutraequ}) converge $\mathbb{P}$-a.s. and
in $L^{p}_{\mathbb{P}}$, as $n\to\infty$, to the orthogonal projection
of $X$ on the (in the ergodic case, at most one-dimensional) space of functions
$Y\in L_{\mathbb{P}}^{p}$ with $TY=e^{-i\theta}Y$.}
\qed

\medskip

The following lemma is a corollary of the previous one. Its proof is
basically the same as that of the equality (16) in
\cite{CuMePe}.\hskip.2pt\footnote{There is a typo in \cite{CuMePe}:
according to the notation there (avoid confusion with our notation) the
correct statement is the following: if $e^{-2it}$ is not an eigenvalue
of $\theta$, (16) is valid $\mathbb{P}$-a.s. (not
$\tilde{\mathbb{P}}$-a.s.) for \textit{every} fixed $u$. An analysis of the
proof shows that the convergence is valid also in the $L^{1}_{
\mathbb{P}}$-sense, which is not explicitly stated there.}

\medskip

\begin{lemmaapp}
\label{rel16cunmerpel}
With the notations and definitions given in Section~\ref{genset} and
Section~\ref{resandcom}, assume that $\theta\in[0,2\pi)$ is such that
$e^{-2i\theta}\notin \mathrm{Spec}_{p}(T)$ and let $Y_{0}\in L^{2}_{
\mathbb{P}}$ be given. Then for every $\mathbf{x}=(x_{1},x_{2})$
\[
\lim_{n}\frac{1}{n}\sum_{k=0}^{n-1}
E_{k-1}\bigl(\mathbf{x}\cdot\bigl(T^{k}Y _{0}e^{ik\theta}
\bigr)\bigr)^{2}= \frac{(x_{1}^{2}+x_{2}^{2})}{2}E\llvert Y_{0}
\rrvert^{2}, \qquad\mathbb{P}\mbox{-a.s. and in }L^{1}_{\mathbb{P}},
\]
where the (probability one) set of pointwise convergence does not depend
on $x_{1}$, $x_{2}$.
\end{lemmaapp}

{\bf Proof:} Notice that
\[
E_{k-1}\bigl(\mathbf{x}\cdot\bigl(T^{k}Y_{0}e^{ik\theta}
\bigr)\bigr)^{2}=T^{k}E_{-1}\bigl( \mathbf{x}\cdot
\bigl(Y_{0}e^{ik\theta}\bigr)\bigr)^{2}
\]
and adapt the argument leading to (16) in \cite{CuMePe} (alternatively, see the proof of Lemma~5 in \cite{bardis}).
\qed

\bigskip 

\begin{remarkapp}
\label{remsetofcon}
{The set of probability one in Lemma~\ref{rel16cunmerpel} can be
described as the set where the pointwise convergences
 $$\frac{1}{n}\sum_{k=0}^{n-1}T^{k}E_{-1}
\bigl[\bigl(\mbox{\upshape Re}(Y_{0})\bigr)^{2}\bigr]e^{i2k
\theta}\to_{n} 0, \quad\quad   \frac{1}{n}\sum_{k=0}^{n-1}T^{k}E_{-1}
\bigl[\bigl(\mbox{\upshape Im}(Y_{0})\bigr)^{2}
\bigr]e^{i2k
\theta}\to_{n} 0,$$ 
$$\frac{1}{n}\sum_{k=0}^{n-1}T^{k}E_{-1}
\bigl[\mbox{\upshape Re}(Y_{0})\mbox{\upshape Im}(Y_{0})
\bigr]e^{i2k
\theta}\to_{n} 0,  \quad\quad \frac{1}{n}\sum
_{k=0}^{n-1}T^{k}E_{-1}\llvert
Y_{0}\rrvert^{2} \to_{n} E\llvert Y_{0}
\rrvert^{2}$$
hold (for fixed versions of the functions involved).} The details are
left to the reader (alternatively, see the proof of Lemma~5 in
\cite{bardis}).
\end{remarkapp}

 \bigskip

The following lemma has a very classical flavor but it is not visible in the literature. We use it in Section \ref{genset} to understand systematically the nature of the quenched results obtained in the paper, and for the proof of Corollary \ref{marcasdou}.

\medskip

\begin{lemmaapp}
\label{mealem}
Let $(\Theta,\mathcal{B}, \lambda)$ and $(\Omega,\mathcal{F},
\mathbb{P})$ be probability spaces and let $\mathcal{F}_{0}\subset
\mathcal{F}$ be a sigma algebra such that $E[\cdot|\mathcal{F}
_{0}]$ admits a regular version in the sense explained in Section~\ref{genset}: there exists a family of probability measures
$\{\mathbb{P}_{\omega}\}_{\omega\in\Omega}$ such that for every
(version of) $Y\in L^{1}_{\mathbb{P}}$,
\[
\omega\mapsto
\int_{\Omega}Y(z)\,d\mathbb{P}_{\omega}(z)
\]
defines an $\mathcal{F}_{0}$-measurable version of $E[Y|\mathcal{F}
_{0}]$. Define, for any (version of) $f\in L^{1}_{\lambda\times
\mathbb{P}}$,
\[
\tilde{f}(\theta,\omega):=
\int_{\Omega}f(\theta,z)\,d\mathbb{P}_{
\omega}(z),
\]
provided that the integral exists, and zero otherwise. Then
$\tilde{f}$ is a version of $E[f|\mathcal{B}\otimes\mathcal{F}_{0}]$.
\end{lemmaapp}

\medskip

{\bf Proof:} Note that, by Fubini's theorem, $f(\theta,\cdot)\in L
^{1}_{\mathbb{P}}$ for $\lambda$-a.e. $\theta$, and that for such
$\theta$, $\tilde{f}(\theta,\cdot)$ is a version of $E_{0}[f(
\theta,\cdot)]$. If we succeed proving that $\tilde{f}$ is
$\mathcal{B}\otimes\mathcal{F}_{0}$-measurable it follows by Fubini's
theorem again that given any rectangular set $E=A\times B\in
\mathcal{B}\otimes\mathcal{F}_{0}$
\[
\int_{E}f(\theta,\omega)\,d(\lambda\times\mathbb{P})=
\int_{A}
\int _{B}f(\theta,\omega)\,d\mathbb{P}(\omega)\,d\lambda(
\theta)=
\int _{A}
\int_{B}\tilde{f}(\theta,\omega)\,d\mathbb{P}(\omega)\,d
\lambda (\theta),
\]
so that, by an application of the $\pi-\lambda$ theorem similar to the
one at the end of this proof,

\[
\int_{E}{f}(\theta,\omega)\,d(\lambda\times\mathbb{P})=
\int_{E} \tilde{f}(\theta,\omega)\,d(\lambda\times\mathbb{P})
\]
for all $E\in\mathcal{B}\otimes\mathcal{F}_{0}$. Thus it suffices to
prove the $\mathcal{B}\otimes\mathcal{F}_{0}$-measurability of
$\tilde{f}$.

Note that if $f^{+}:=fI_{[f\geq0]}$ and $f^{-}:=-fI
_{[f<0]}$ are, respectively, the nonnegative and negative parts of a
real-valued $f\in L^{1}_{\lambda\times\mathbb{P}}$, we can recover the
definition of $\tilde{f}$ via the formula
\[
\tilde{f}=\bigl(\tilde{f^{+}}-\tilde{f^{-}}
\bigr)I_{[\widetilde{|f|}>0]}
\]
(here $\widetilde{|f|}:=\tilde{g}$ with $g:=|f|$), and we can create a (measurable) formula for an arbitrary $f\in L^{1}
_{\lambda\times\mathbb{P}}$ by applying this to its real and imaginary
parts.

Thus it suffices to assume that $f$ is nonnegative. We will do so for
the rest of the proof.

It is well known that every nonnegative function $f$ can be approximated
by simple functions $f_{n}$ with $f_{n}$ increasing to $f$
(\cite{BilProMea}, Theorem~13.5, p.~185). Then, by the monotone
convergence theorem,
\[
\tilde{f}=\lim_{n}(\tilde{f}_{n}I_{A_{f}})
\]
where
\[
A_{f}:=\bigcup_{n\in\mathbb{N}}\bigcap
_{k\in\mathbb{N}}[\tilde{f} _{k}\leq n] .
\]
Thus, it suffices to see that for each simple function $f$,
$\tilde{f}$ is $\mathcal{B}\otimes\mathcal{F}_{0}$-measurable, and
therefore it suffices (by linearity) to prove this if $f=I_{E}$ for any
$E\in\mathcal{B}\otimes\mathcal{F}$.

Let us do it: if $E=A\times B$ is a rectangular set, then
\[
\tilde{I}_{E}(\theta,\omega)=I_{A}(\theta)
\mathbb{P}_{\omega}(B) ,
\]
which is clearly $\mathcal{B}\otimes\mathcal{F}_{0}$ measurable. It
follows that $\tilde{I}_{E}$ is
$\mathcal{B}\otimes\mathcal{F}_{0}$-measurable if $E$ is any finite
union of disjoint rectangles in $\mathcal{B}\otimes\mathcal{F}$.

Now consider the family $\tilde{\mathcal{F}}'$ of sets $E \in
\mathcal{B}\otimes\mathcal{F}$ such that $\tilde{I}_{E}$ is
$\mathcal{B}\otimes\mathcal{F}_{0}$-measurable. Since for any family
$\{E_{n}\}_{n}{\subset}\tilde{\mathcal{F}}'$ of mutually disjoint sets
\[
\tilde{I}_{\bigcup_{n}E_{n}}=\sum_{n}
\tilde{I}_{E_{n}}
\]
(apply the monotone convergence theorem) and $\Theta\times\Omega$
 is a subset of
$\tilde{\mathcal{F}}'$,
$\tilde{\mathcal{F}}'$ is a
$\lambda$-system. Since $\tilde{\mathcal{F}}'$ includes the finite
unions of disjoint rectangles it follows, by the $\pi-\lambda$
theorem, that $\tilde{\mathcal{F}}'=\mathcal{B}\otimes\mathcal{F}$. \qed


\medskip

\begin{corapp}
\label{cormealem}
Under the conditions of Lemma \textup{\ref{mealem}}, if $\mathcal{B}_{0}\subset
\mathcal{B}$ is a sigma algebra such that $E[\,\cdot\,|\mathcal{B}
_{0}]$ is regular with regular measures $\{\lambda_{\theta}\}_{
\theta\in\Theta}$, so that
\[
\theta\mapsto
\int_{\Theta} g(z) \,d\lambda_{\theta}(z)
\]
is $\mathcal{B}_{0}$-measurable and defines a version of $E[\,g|
\mathcal{B}_{0}]$ for any $g\in L^{1}_{\lambda}$, then $\{
\lambda_{\theta}\times\mathbb{P}_{\omega}\}_{(\theta,\omega)
\in\Theta\times\Omega}$ is a family of regular measures for the conditional
expectation $E[\,\cdot\,|\mathcal{B}_{0}\otimes\mathcal{F}_{0}]$ with respect to $\mathcal{B}_{0}\otimes\mathcal{F}_{0}$ and
$\lambda\times\mathbb{P}$.
\end{corapp}

{\bf Proof:} Given any function $f=f(\theta,\omega)\in L^{1}_{
\lambda\times\mathbb{P}}$\vspace*{-2pt}
\[
E[f|\mathcal{B}_{0}\otimes\mathcal{F}_{0}]=E\bigl[E[f|
\mathcal{B} \otimes\mathcal{F}_{0}]|\mathcal{B}_{0}\otimes
\mathcal{F}_{0}\bigr]=E[ \tilde{f}|\mathcal{B}_{0}\otimes
\mathcal{F}_{0}],
\]
where $\tilde{f}$ is the function specified by Lemma~\ref{mealem}. A~second application of this lemma with $\mathcal{B}\otimes\mathcal{F}
_{0}$ in the role of $\mathcal{B}\otimes\mathcal{F}$ and with
$\mathcal{B}_{0}$ in the role of $\mathcal{F}_{0}$ gives that\vspace*{-2pt}
\[
(\theta,\omega)\mapsto
\int_{\Theta}
\int_{\Omega}f(x,z) \,d \mathbb{P}_{\omega}(z)\,d
\lambda_{\theta}(x)
\]
defines a $\mathcal{B}_{0}\otimes\mathcal{F}_{0}$-measurable version of $E[f|\mathcal{B}_{0}\otimes\mathcal{F}_{0}]$.
\qed

\bigskip
 
We finish this section by stating two widely known limit theorems, which
are building blocks of the results presented in this paper.

\medskip
\begin{thmapp}[The Lindeberg--L\'{e}vy theorem for martingales]
\label{bilmarclt}
For each $n\in\mathbb{N}^{*}$, let $\Delta_{n1},\dots,\Delta_{nk},
\dots$ be a sequence of real-valued martingale differences with respect
to some increasing filtration $\mathcal{F}_{0}^{n}\subset\cdots
\subset\mathcal{F}_{k}^{n}\subset\cdots$.

Define, for every $(n,k)\in\mathbb{N}^{*}\times \mathbb{N}^{*}$, $\sigma_{nk}^{2}:=E[\Delta_{nk}^{2}\|\mathcal{F}^{n}_{k-1}]$, and assume that $\sum_{k}\Delta_{nk}$ and $\sum_{k}\sigma_{nk}^{2}$ converge with probability one. If for some $\sigma\geq0$ the following two conditions hold
\begin{enumerate}
\item[1.] $\sum_{k\geq1}\sigma_{nk}^{2}\Rightarrow_{n}\sigma^{2}$,
\item[2.] $\sum_{k\geq1}E[\Delta_{nk}^{2}I_{[\Delta_{nk}\geq\epsilon]}]
\to_{n}0$, for every $\epsilon>0$,
\end{enumerate}
then $Z_{n}:=\sum_{k\geq1}\Delta_{nk}\Rightarrow_{n}\sigma N$ where
$N$ is a standard normal random variable.
\end{thmapp}

{\bf Proof:}
\cite{BilProMea}, p.~477.

\medskip
\begin{thmapp}[The functional form of Theorem~\ref{bilmarclt}]
\label{bilinvpri}
For each $n\in\mathbb{N}^{*}$, let $\Delta_{n1},\dots,\Delta_{nk},
\ldots $ be a sequence of real-valued martingale differences with respect
to some increasing filtration $\mathcal{F}_{0}^{n}\subset\cdots
\subset\mathcal{F}_{k}^{n}\subset\cdots$.

Define, for every $(n,k)\in\mathbb{N}^{*}\times \mathbb{N}^{*}$, $\sigma_{nk}^{2}:=E[\Delta_{nk}^{2}\|\mathcal{F}^{n}_{k-1}]$. If for some $\sigma\geq0$ the following two conditions hold
for every $t\geq0$, $\epsilon>0$
\begin{enumerate}
\item[1.] $\sum_{k\leq nt}\sigma_{nk}^{2}\Rightarrow_{n} \sigma^{2} t$,
\item[2.] $\sum_{k\leq nt}E[\Delta_{nk}^{2}I_{[\Delta_{nk}\geq\epsilon
]}]\to
_{n} 0$,
\end{enumerate}
then the random functions $X_{n}(t):=\sum_{k\leq nt}\Delta_{nk}$
converge in distribution, as $n\to\infty$, to $\sigma W$ in the sense of $D[[0,\infty)]$
where $W$ is a standard Brownian motion.
\end{thmapp}

{\bf Proof:} This is a slight reformulation of Theorem~18.2 in
\cite{Bilconpromea}, (pp.~194--195): the case $\sigma>0$ follows by a
simple renormalization, and to cover the case $\sigma=0$, note that the
convergence (18.6) in \cite{Bilconpromea} becomes a simple consequence
of the definition given there of $\zeta_{nk}$ and the hypothesis
(corresponding to $\sigma=0$)\vspace*{-2pt}
\[
\sum_{k\leq nt}\sigma_{nk}^{2}
\Rightarrow_{n}0
\]
for every $t\geq0$.\qed 

\subsection*{Acknowledgements}
This work is part of the research project carried out by the author
under the tutelage of M.~Peligrad, and he is indebted to her for
suggesting accurate references and giving important comments, in
particular with regards to the scope of our methods under the weakened
Hannan condition (\ref{imhancon}).

My conversations with C. Dragan were very useful. Inspirational
suggestions were given by professors W. Bryc and Y. Wang during our
meetings at the University of Cincinnati Probability Seminar.

It is also difficult to overestimate the role of the anonymous referee
with regards to improving the presentation and the results presented
along this paper. Besides many stylistic corrections, the author would
like to acknowledge that the referee's suggestions gave rise to the
ideas leading to the proof of Theorem~\ref{barinvprihancon} under
(\ref{imhancon}) presented here\footnote{The previous proof departed
from a more elementary level but it was also longer, the reader can
consult it in \cite{bardis}.} and the corresponding result under
condition (\ref{maxwoo}). The referee also suggested to look at
(\ref{equpromaxwooae}) as a suitable condition for the $\lambda$-a.e.
fulfillment of (\ref{maxwoo}), providing a reference (\cite{ass}) that led to the materialization of Proposition~\ref{proconsqu} and its
consequences.

The author was supported by the NSF Grants DMS-1208237, DMS-1512936, and
by a Laws scholarship for the majority of the time of research and
writing of this paper. The revision following the last report by the
referee was done during the first two months of the author's stay at the
\'{E}cole Polytechnique, supported by a grant from the Risk Foundation.

\end{document}